\documentclass[12pt,twoside]{amsart}\usepackage{amssymb}
\usepackage{enumerate}

\newtheorem{thmspec}{\relax}

\newtheorem{theorem}{Theorem}[section]
\newtheorem{thm}[theorem]{Theorem}
\newtheorem{lem}[theorem]{Lemma}
\newtheorem{cor}[theorem]{Corollary}
\newtheorem{prop}[theorem]{Proposition}

\theoremstyle{definition}

\theoremstyle{remark}

\numberwithin{equation}{section}
\def \Bbb{\mathbb}

\def\onto{{\kern3pt\to\kern-8pt\to\kern3pt}}

\def\<{\langle}
\def\>{\rangle}
\def\|{{\ |\ }}

\def\onto{\twoheadrightarrow}

\def\-{\underline}

\def\env{\operatorname{env}}

\def\dist{\operatorname{dist}}
\def\Re{\operatorname{Re}}

\def\diam{\operatorname{diam}}

\def\N{\Bbb N}

\def\R{\Bbb R}
\def\C{\Bbb C}

\def\X{\Bbb X}

\def\<{\langle}
\def\>{\rangle}

\catcode`\@=11
\def\serieslogo@{\relax}
\def\@setcopyright{\relax}
\catcode`\@=12

\title[ A boundary cross theorem ]
{  A boundary cross theorem\\
 for separately holomorphic  functions
}

\begin{document}

\author{Peter Pflug}
\address{Peter Pflug\\
Carl von  Ossietzky Universit\"{a}t Oldenburg \\
Fachbereich  Mathematik\\
Postfach 2503, D--26111\\
 Oldenburg, Germany}
\email{pflug@mathematik.uni-oldenburg.de}

\author{   Vi\d{\^e}t-Anh  Nguy\^en}
\address{  Vi\d{\^e}t-Anh  Nguy\^en \\
Carl von  Ossietzky Universit\"{a}t Oldenburg \\
Fachbereich  Mathematik\\
Postfach 2503, D--26111\\
 Oldenburg, Germany}
\email{nguyen@mathematik.uni-oldenburg.de}

\subjclass[2000]{Primary 32D15, 32D10}
\date{}

\keywords{ Cross Theorem, holomorphic extension, plurisubharmonic measure.}

\begin{abstract}
  Let $D\subset \C^n,$ $G\subset \C^m$ be pseudoconvex domains, let
  $A$ (resp. $B$) be an  open subset of the boundary $\partial D$ (resp.
  $\partial G$) and let $X$ be the  $2$-fold cross $((D\cup A)\times B)\cup (A\times(B\cup G)).$
  Suppose in addition  that  the domain $D$ (resp. $G$) is  {\it locally  $\mathcal{C}^2$ smooth
  on $A$} (resp. $B$).
   We shall determine the "envelope of holomorphy"
  $\widehat{X}$ of $X$ in the sense that any function continuous on $X$ and separately
  holomorphic
  on $(A\times G) \cup (D\times B)$ extends to a function   continuous on $\widehat{X}$
  and holomorphic  on the interior of $\widehat{X}.$ A generalization of this result
  for an $N$-fold cross is also given.
\end{abstract}
\maketitle

\section{ Introduction and statement of the main results }

In order to recall here the classical cross theorem and to state our
results, we need to introduce some notation and terminology.
 In fact, we keep the main notation from the book by  Jarnicki and the
 first author \cite{jp1} and from the survey article by Sadullaev
\cite{sa1}.
\subsection{Plurisubharmonic measures}
Let $\Omega\subset \C^n$ be an open set.  For any function $u$ defined on $\Omega,$ let
\begin{equation*}
 \hat{u}(z):=
\begin{cases}
u(z),
  & z\in   \Omega,\\
 \limsup\limits_{w\in \Omega,\ w\to z}u(w), & z \in \partial \Omega.
\end{cases}
\end{equation*}
For a set  $A\subset \overline{\Omega}$ put
\begin{equation*}
h_{A,\Omega}:=\sup\left\lbrace u\ :\  u\in\mathcal{PSH}(\Omega),\ u\leq 1\ \text{on}\ \Omega,\
   \hat{u}\leq 0\ \text{on}\ A    \right\rbrace,
\end{equation*}
where $\mathcal{PSH}(\Omega)$ denotes the set   of all functions  plurisubharmonic
on $\Omega.$

 We first suppose that $\Omega$ is bounded. Then
 the {\it plurisubharmonic measure} of $A$ relative to $\Omega$ is given
 by
\begin{equation}
 \omega(z,A,\Omega):=
 \widehat{h^{\ast}_{A,\Omega}}(z)    , \qquad z \in  \Omega\cup A,
\end{equation}
where $u^{\ast}$ denotes the  upper semicontinuous regularization of a function $u.$

From now on let $\Omega$ be an arbitrary (not necessarily bounded) open set and
we shall be concerned with the following two cases.

\smallskip

\noindent{\bf Case I:} $A\subset \Omega.$

In this case, we define the {\it plurisubharmonic measure} of $A$ relative to $\Omega$
 as follows.
\begin{equation*}
\omega(\cdot,A,\Omega):=\lim\limits_{k\to +\infty} h^{\ast}_{A\cap
\Omega_k,\Omega_k},
\end{equation*}
where $\left ( \Omega_k \right )^{\infty}_{k=1}$ is a sequence
of relatively compact open sets $\Omega_k\subset \Omega_{k+1}\Subset
\Omega$ with $\cup_{k=1}^{\infty} \Omega_k=\Omega.$
Observe that the definition is independent  of the exhausting sequence
$\left ( \Omega_k \right )^{\infty}_{k=1}.$  Moreover,
$\omega(\cdot,A,\Omega) \in \mathcal{PSH}(\Omega).$

\smallskip

\noindent{\bf Case II:} $A$ is an open subset of $ \partial\Omega.$

We suppose in addition that
 $\Omega$ is  {\it locally $\mathcal{C}^2$ smooth} on $A$ (i.e.
for any $\zeta\in A,$ there exist an open neighborhood $U=U_{\zeta}$ of
 $\zeta$ in $\C^n$ and a   real function
 $\rho=\rho_{\zeta}\in \mathcal{C}^2(U)$ such that $\Omega\cap U=\lbrace z\in U:\
 \rho(z)<0\rbrace$ and $d\rho(\zeta)\not=0$).

In this case,  the {\it plurisubharmonic measure} of $A$ relative to $\Omega$
is  a function on $\Omega\cup A$ given by
\begin{equation*}
\omega(z,A,\Omega):=\lim\limits_{k\to +\infty} \omega\left(z, A_k,\Omega_k\right),
\qquad z\in \Omega\cup A,
\end{equation*}
where the function  $  \omega\left(\cdot, A_k,\Omega_k\right)$ is given by (1.1) and
 $\left ( \Omega_k \right )^{\infty}_{k=1}$ is a sequence
of relatively compact open sets $\Omega_k \Subset
\C^n$  and $\left ( A_k\right )^{\infty}_{k=1}$ is a sequence of open subsets of $A$  such that
\begin{itemize}
\item[(i)] $\Omega_k\subset \Omega_{k+1}$ and $\cup_{k=1}^{\infty} \Omega_k=\Omega;$
\item[(ii)]  $A_k\subset A_{k+1}$ and $A_k\subset \partial \Omega\cap\partial \Omega_k$
and $\cup_{k=1}^{\infty} A_k=A;$
 \item[(iii)]  for any point $\zeta\in A$ there is an open neighborhood
 $V=V_{\zeta}$ of $\zeta$ in $\C^n$ such that $V\cap \Omega$ =$V\cap\Omega_k$ for
 some $k.$
  \end{itemize}
  In Section 3 below, we
 shall prove
 that the definition is independent  of the chosen exhausting sequences
$\left ( \Omega_k \right )^{\infty}_{k=1}$ and  $\left ( A_k \right )^{\infty}_{k=1}.$ Moreover,
$\omega(\cdot,A,\Omega)|_{\Omega} \in \mathcal{PSH}(\Omega).$

\subsection{Cross and separate holomorphicity}
 Let $N\in\N,\ N\geq 2,$ and let $\varnothing\not = A_j\subset
\overline{D_j}\subset\C^{n_j},$ where $D_j$ is a domain, $j=1,\ldots,N.$ We define
an {\it $N$-fold cross} $X,$ {\it its interior} $X^{\text{o}}$ and a new set $A$ as
\begin{eqnarray*}
X &=&\X(A_1,\ldots,A_N; D_1,\ldots,D_N)\\
&:=& \bigcup_{j=1}^N  A_1\times\cdots\times A_{j-1}\times (D_j\cup A_j)\times
A_{j+1}\times\cdots \times A_N\subset \C^{n_1+\cdots+n_N} =\C^n,\\
X^{\text{o}} &=&\X^{\text{o}}(A_1,\ldots,A_N; D_1,\ldots,D_N)\\
&:=& \bigcup_{j=1}^N  A_1\times\cdots\times A_{j-1}\times D_j\times
A_{j+1}\times\cdots \times A_N,\\
A&:=&A_1\times \cdots\times A_N.
\end{eqnarray*}
In particular, if $A_j\subset D_j,\ j=1,\ldots,N,$ then we have $X=X^{\text{o}}.$
Moreover put
\begin{equation*}
\omega(z):=\sum\limits_{j=1}^N \omega(z_j,A_j,D_j),\qquad
z=\left(z_1,\ldots,z_N\right)\in (D_1\cup A_1)\times
\cdots\times (D_N\cup A_N).
\end{equation*}

For an $N$-fold cross $X:=\X(A_1,\ldots,A_N;D_1,\ldots,D_N)$
 let
\begin{equation*}
\widehat{X}:=\left\lbrace z=\left(z_1,\ldots,z_N\right )\in (D_1\cup A_1)\times
\cdots\times (D_N\cup A_N):\ \omega(z)  <1
\right\rbrace.
\end{equation*}
Then the set of all interior points  of
$\widehat{X}$ is given by
\begin{equation*}
 \widehat{X}^{\text{o}} :=\left\lbrace z= \left(z_1,\ldots,z_N\right )\in D_1\times
\cdots\times D_N:\  \omega(z)<1
\right\rbrace.
\end{equation*}
We say that a function $f:\ X\longrightarrow \C$ is  {\it separately holomorphic}
on $X^{\text{o}} $ and
write $f\in\mathcal{O}_s(X^{\text{o}}),$   if
 for any $j\in\lbrace
1,\ldots,N\rbrace$ and $(a^{'},a^{''})\in  (A_1\times\cdots\times
 A_{j-1})\times(A_{j+1}\times \cdots \times A_N) $
 the function $f(a^{'},\cdot,a^{''})|_{D_j}$ is holomorphic  on $D_j.$

 Finally, throughout the paper,
the notation  $\vert f\vert_M$ denotes $\sup_M \vert f\vert.$

\subsection{Motivations for our work}
We are now able to formulate what we will quote in the sequel as the {\it classical cross
theorem}.
\renewcommand{\thethmspec}{Theorem 1}
  \begin{thmspec} (Alehyane \& Zeriahi \cite{az})
Let  $D_j\subset\C^{n_j}$ be a pseudoconvex domain and  $ A_j\subset
D_j$ a locally  pluriregular subset,
  $j=1,\ldots,N.$
Then for any function $f\in\mathcal{O}_s(X),$ there is a unique function
$\hat{f}\in\mathcal{O}(\widehat{X})$
such that $\hat{f}=f$ on $X.$
\end{thmspec}
There is a long list  of papers dealing with this theorem under various
assumptions. For a historical discussion, see the survey article \cite{pf1}.

The question naturally arises  how  the situation changes when the sets $A_j$ live on the
boundary $\partial D_j,$  $j=1,\ldots,N.$

The first results in this direction are obtained by Malgrange--Zerner \cite{ze},
Komatsu \cite{ko} and Dru\.{z}kowski \cite{dr}, but only for  some special crosses.
Recently, Gonchar \cite{go1,go2} has proved  the following remarkable more
general theorem.

\renewcommand{\thethmspec}{Theorem 2}
  \begin{thmspec}
Let  $D_j\subset\C$ be a   Jordan domain and let  $\varnothing\not= A_j$  be an open set of the boundary
$\partial D_j,$   $j=1,\ldots,N.$
Then for any  function $f\in\mathcal{C}(X)\cap\mathcal{O}_s(X^{\text{o}})$ there is a unique function
$\hat{f}\in\mathcal{C}(\widehat{X})\cap\mathcal{O}(\widehat{X}^{\text{o}})$
such that $\hat{f}=f$ on $X.$ Moreover, if $\vert f\vert_X<\infty$ then
\begin{equation}
 \vert \hat{f}(z)\vert\leq \vert f\vert_A^{1-\omega(z)} \vert
 f\vert_X^{\omega(z)},\qquad z\in\widehat{X}.
\end{equation}
\end{thmspec}
It should be observed that under the hypothesis of Theorem 2 one has $X\subset\widehat{X}.$
We remark that the original formulation of Gonchar is slightly different from Theorem 2.
However, his proof is still valid also in this new formulation.

The main purpose of this work is to generalize Gonchar's theorem to higher
dimensions.
\subsection{The main result}
We are now ready to state the main  result.

\renewcommand{\thethmspec}{Main Theorem}
  \begin{thmspec}
Let  $D_j\subset\C^{n_j}$ be a pseudoconvex domain and let  $\varnothing\not = A_j$
be an open set of $ \partial D_j,$  $j=1,\ldots,N.$ Suppose in addition
that each domain $D_j$ is  locally
$\mathcal{C}^2$ smooth on $A_j,$   $j=1,\ldots,N.$
 Then  $X\subset \widehat{X}$ and for any function $f\in  \mathcal{C}(X)\cap\mathcal{O}_s(X^{\text{o}})       ,$ there is a unique function
$\hat{f}\in\mathcal{C}(\widehat{X})\cap\mathcal{O}(\widehat{X}^{\text{o}}) $
such that $\hat{f}=f$ on $X.$
 Moreover, if $\vert f\vert_X<\infty$ then
 \begin{equation}
 \vert \hat{f}(z)\vert \leq \vert f\vert_A^{1-\omega(z)} \vert
 f\vert_X^{\omega(z)},\qquad z\in\widehat{X}.
\end{equation}
\end{thmspec}

We will give here a short outline of the proof.

The main idea is to combine Gonchar's theorem and the classical cross theorem
with the {\it slicing method}. More precisely,
for each domain $D_j$ we shall associate a
family of     $\mathcal{C}^2$ smooth planar domains which are, roughly speaking,
the intersection of  an open tubular  neighborhood of $A_j$ in $D_j\cup A_j$ with the family of  normal complex lines to $A_j$ parameterized by $A_j.$
 One important property
of this family is that the harmonic measures for its domains  depend, in
some sense, continuously on the parameter of  $A_j.$
Another important property is that    there is a relation between
the plurisubharmonic measure of $D_j$ and the harmonic
measure of the domains in the above family.
 Applying Gonchar's theorem and the slicing method, we
shall
find an extension $\hat{f}$
such that $\hat{f}$ is holomorphic on a subdomain of each domain in this family.
The two important properties mentioned above, combined with a variant of the classical cross theorem,
 will allow us to propagate the holomorphicity  from these one-dimensional subdomains
  to the desired envelope of holomorphy.

This paper is organized as follows.

We begin Section 2 by
collecting some background of the potential   theory and some classical results.
Next we establish a uniform  estimate for the Poisson kernels
which will play an important role in the proof of the main theorem.

Based on the results of Section 2,
Section 3 develops necessary estimates for the plurisubharmonic measures
that will be used later in Section 5.

 Section 4 provides the first step of the proof.
More precisely, on the one hand
we will consider the following mixed situation where there is at least one
index $j$ such that the
factor $A_j$ of the cross $X$ is inside $D_j.$
On the other hand, we will establish  some quantitative versions
of the classical cross theorem.

Section 5 establishes  the main theorem in the case of a $2$-fold cross.

The complete proof of the main theorem will be given in Section 6 together with
some concluding remarks and open questions.

\smallskip

\indent{\it{\bf Acknowledgment.}} The paper was written while  the second author was visiting the
Carl von  Ossietzky Universit\"{a}t Oldenburg
 being  supported by The Alexander von Humboldt
Foundation. He wishes to express his gratitude to these organizations.

\section{Auxiliary results}
\subsection{Harmonic measure for a planar domain}

We recall   some classical facts from the book of Ransford \cite{ra}.
Let $D$ be a proper subdomain of $\C\cup\{\infty\}$ such that the boundary $\partial D$
(with respect to $\C\cup\{\infty\}$)
is non-polar. Let $\mathcal{P}_D$ be the Poisson projection of $D$
and $A$ a Borel subset of $\partial D.$ Consider the bounded function
\begin{equation*}
 1_{\partial D\setminus A}(\zeta):=
\begin{cases}
1,
  &\zeta\in  \partial D\setminus A,\\
 0, & \zeta \in A.
\end{cases}
\end{equation*}
 Then, by Theorem 4.3.3 of \cite{ra}, the {\it  harmonic measure} of the set
 $\partial D\setminus A$ (or equivalently $h_{A,D}$) is exactly the
 Perron solution of the generalized Dirichlet problem with boundary data $1_{\partial D\setminus
 A}.$ In other words, one has
\begin{equation}
 h_{A,D}=\mathcal{P}_D[1_{\partial D\setminus
 A}].
\end{equation}
\subsection{A uniqueness theorem and a Two-Constant Theorem}

The following  uniqueness theorem is very useful.

\begin{thm} (see \cite{hc})  Let $D\subset \C^n$ be a domain such that $D$ is locally  $\mathcal{C}^1$ smooth on some
open set $U$ of   $\partial D.$
If a set $E\subset D\cup U$ has  a positive $(2n-1)$-dimensional
Hausdorff measure, then, for $f\in\mathcal{C}(D\cup U)\cap\mathcal{O}(D),$
$f=0$ on $E$ implies $f\equiv 0.$
\end{thm}
\begin{proof} The only nontrivial case is that $E\subset U.$
In this case by taking the intersection of $D$ with a bundle of complex
lines and applying the  classical one-dimensional boundary uniqueness theorem of Privalov,
one may find a set $E^{'}\subset D$ close to $E$ such that
$E^{'}$ has  a positive $2n$-dimensional
Hausdorff measure and  $f=0$ on $E^{'}$. This completes the proof.
\end{proof}
The following Two-Constant Theorem for plurisubharmonic  functions
will play a vital role in this paper.
\begin{thm}
If $u$ is a  plurisubharmonic function in a bounded open set $D\subset\C^n$ and
$u\leq M$ on $D$ and $\hat{u}\leq m$ on some subset $A$ of $\overline{D},$
then
\begin{equation*}
\hat{u}(z)\leq m(1-\omega(z,A,D))+M\cdot\omega(z,A,D),\qquad z\in \overline{D}.
\end{equation*}
\end{thm}
\begin{proof}
It follows immediately from the definition of $\omega(\cdot,A,D)$ given in
Subsection 1.1.
\end{proof}

\subsection{Uniform estimate for the Poisson kernels of a family of
 $\mathcal{C}^2$ smooth domains}
 In what follows we fix an integer $N\geq 2$ and let
 $\dist(\cdot,\cdot)$ denote the Euclidean distance and let
 $B(a,r)$ $(a\in\R^N,r>0)$ denote the Euclidean ball of center $a$ and radius $r.$
 We say that a domain  $ D\subset \R^N$ is  {\it $\mathcal{C}^2$ smooth}
 if $D$ is bounded and  admits a defining function $\rho \in\mathcal{C}^2(\R^N)$
 such that $d\rho(z)\not=0$ for all $z\in\partial D.$
 Let $P_D$ denote its Poisson kernel.
  We begin this subsection with the following result
due to N. Kerzman (see \cite{ke} and \cite{kr}).
\begin{thm}
Let $D\subset \R^N$ be a $\mathcal{C}^2$ smooth domain   that satisfies
\begin{equation*}
\diam(D):=\sup\limits_{x,y\in D} \vert x-y\vert \leq M \quad\text{for
some finite constant}\ M.
\end{equation*}
Then

 \noindent 1) there is a  positive number
$r=r(D) $ such that
 for each $y\in\partial
D$ there are balls $B(c_y,r)\subset D$
and  $B(\widetilde{c}_y, r )\subset\R^N\setminus \overline{D}$
that satisfy
\begin{eqnarray*}
\overline{B(\widetilde{c}_y,r)}\cap
\overline{D} &=&\{y\},\\
\overline{B(c_y,r)}\cap (\R^N\setminus D)&=&\{y\};
\end{eqnarray*}

\noindent 2) there is a finite constant $C$ which depends only on $N, r,$ and
$M$ such that
\begin{equation*}
P_{D}(x,y)\leq C\frac{\dist (x,\partial D)}
{\vert x-y\vert^N},\qquad x\in D,\ y\in \partial D.
\end{equation*}
\end{thm}
\begin{proof} This theorem  is implicitly  proved in Lemmas 8.2.3--8.2.5 and Proposition 8.2.6
in the book of Krantz \cite{kr}. We only mention here that
Kerzman's idea is to compare the Green function and the Poisson kernel for $D$
with the corresponding functions  for the internally and externally
tangent balls $B(c_y,r)$
and  $B(\widetilde{c}_y,r) $ (and also for their complement).
\end{proof}

Now we reformulate Kerzman's theorem in order to obtain
an uniform upper bound for the Poisson kernels of a family of domains which
depend, in some sense, continuously on a parameter.

\begin{cor} Let
 $\left(D_{\alpha}\right)_{\alpha\in I}$
be a   family of $\mathcal{C}^2$ smooth domains in $\R^N$ indexed by a set $I.$
Suppose that

\noindent 1) there is a finite constant $M$   such that
for all $\alpha\in I,$
\begin{equation*}
\diam(D_{\alpha})  \leq M ;
\end{equation*}

 \noindent 2) there is a finite  positive number
$r $ such that
 for each $\alpha\in I,$ $y\in\partial
D_{\alpha},$ there are balls $B(c_{y,\alpha},r)\subset D_{\alpha}$
and  $B(\widetilde{c}_{y,\alpha}, r )\subset\R^N\setminus \overline{D}_{\alpha}$
that satisfy
\begin{eqnarray*}
\overline{B(\widetilde{c}_{y,\alpha},r)}\cap
\overline{D}_{\alpha} &=&\{y\},\\
\overline{B(c_{y,\alpha},r)}\cap (\R^N\setminus D_{\alpha})&=&\{y\}.
\end{eqnarray*}

Then there exists a finite constant $C$ such
that
\begin{equation*}
P_{D_{\alpha}}(x,y)\leq C\frac{\dist (x,\partial D_{\alpha})}
{\vert x-y\vert^N},\qquad x\in D_{\alpha},\ y\in \partial D_{\alpha},\
\alpha\in I.
\end{equation*}
\end{cor}
\begin{proof}
It follows immediately from Theorem 2.3.
\end{proof}

 We conclude this section with an example of a family of
 $\mathcal{C}^2$ smooth domains satisfying the hypothesis of Corollary 2.4.

  Let $D$ be domain in $\C^n$ such that $  D$ is
 locally $\mathcal{C}^2$ smooth on an open neighborhood of a point
 $P\in\partial D.$ Let $T^{\C}_P$ (resp. $T^{\R}_P$) denote the complex (resp. real)
 tangent hyperplane to
 $\partial D$ at $P$ and $\pi$ (resp. $\pi^{\C}$) the orthogonal projection from $\C^n$ onto
  $T^{\R}_P$ (resp. $T^{\C}_P$).

By an affine  transformation, we may suppose without loss of
generality that $P= 0,$ $T^{\C}_P=\left\lbrace z_1=0\right\rbrace$
and $T^{\R}_P=\left\lbrace \Re z_1=0\right\rbrace.$
Moreover, there are an open neighborhood $U$ of the origin and a function
$\rho\in\mathcal{C}^2(U)$ such that
\begin{equation}
\rho(0)=0, d \rho(0)=(1,0,\ldots,0)\quad\text{and}\quad  U\cap D=
\left\lbrace \rho <0\right\rbrace.
\end{equation}

 For any  domain $V\subset U,$ any  $Q:=(0,z^{'})=(0,z_2,\ldots,z_n)\in
 T^{\C}_P$,
 consider the planar domain
 \begin{equation}
 V_{Q}  :=\env\left(V\cap \left\lbrace  (t,z^{'}),\  t\in \C
 \right\rbrace\right),
 \end{equation}
 where $\env(G)$ denotes the smallest simply connected domain containing
  (a given planar domain) $G,$ in other words, $\env(G)$ is  obtained from $G$ by adding all its
 holes.
 \begin{prop} Under the above hypothesis and notation,
  there are  open neighborhoods $U_1$ of $P$ in $T^{\C}_P,$
    $ U_2$ of $P$ in $T^{\R}_P $ and $U_3$ of $P$ in $\C^n$ and a
 $\mathcal{C}^2$ smooth subdomain $V\subset D$ such that
 \begin{itemize}
 \item[1)]  $U_1=U_2\cap T^{\C}_P;$
 \item [2)]  $\partial V\cap \partial D$ is an open neighborhood
 of $P$ in $\partial D$ and in $\partial V$ and $\pi$
 is one-to-one from $\partial V\cap \partial D$ onto an open neighborhood of $U_2;$
 \item  [3)]$\left( V_{Q} \right)_{Q\in U_1}$ is a   family of
 $\mathcal{C}^2$ smooth planar simply connected domains which satisfies 1) and 2)
 in Corollary 2.4 and $V_Q\subset D;$
 \item  [4)] there is a finite constant $C$ such that for all $Q\in U_1,$
 $z\in V_Q\cap U_3$ and $\zeta\in \partial V\cap\partial D$ satisfying
 $\pi(\zeta)=\pi(z),$
 \begin{equation*}
\dist(z,\partial V_Q)\leq C\cdot \dist(z,\partial D)
\qquad\text{and}\qquad \dist(z,\zeta)\leq C\cdot \dist(z,\partial D);
 \end{equation*}
 in other words, the quantities $\dist(z,\partial V_Q),$  $\dist(z,\zeta)$ and $
 \dist(z,\partial D)$
 are equivalent.
\end{itemize}
\end{prop}

\begin{proof}
Since  $  D$ is
 locally $\mathcal{C}^2$ smooth on an open neighborhood of a point
 $P\in\partial D,$ a geometric argument (see \cite[p. 325]{kr}) shows that
 there is an $r>0$ such that the sphere $\partial B$ is internally tangent
 to $D$ at $P,$ where the ball $B$ is given by $B:=B((-r,0,\ldots,0),r).$

 Consider the following defining function for the ball $B$
 \begin{equation}
\phi(z):=\frac{(x_1+r)^2+y_1^2 +\vert z^{'}\vert^2 -r^2}{2r},\qquad
z=(x_1+iy_1,z^{'})\in\C^n.
 \end{equation}
A straightforward computation gives that $\vert d\phi\vert =1$ on
$\partial B.$ Next fix a  radial function $\psi\in\mathcal{C}_{0}(\C^n)
$ such that $0\leq \psi\leq 1,$ $\psi(z)=1$ for $\vert z\vert\leq 1$ and
$\psi(z)=0$ for $\vert z\vert\geq 2.$

Since $d\rho(0)=d\phi(0),$ we may choose a sufficiently small $\epsilon_0$
such that $0<\epsilon_0<\frac{r}{4}$ and
\begin{equation}
\left\vert (d\rho-d\phi)(z)\right\vert<\frac{1}{8},\qquad
\vert z\vert\leq 2\epsilon_0.
\end{equation}

Now define  for any $0<\epsilon<\epsilon_0,$
\begin{equation}
\begin{split}
\psi_{\epsilon}(z)&:=\psi\left(\frac{z}{\epsilon}\right),\\
\rho_{\epsilon}&:=\phi+\psi_{\epsilon}(\rho-\phi).
\end{split}
\end{equation}

Observe that $\rho_{\epsilon}(z)=\rho(z)$ for $\vert z\vert \leq \epsilon$
and $\rho_{\epsilon}(z)=\phi(z)$ for $\vert z\vert \geq 2\epsilon.$
Moreover using (2.5)--(2.6) and
the identities  $\rho(0)=\phi(0),$ $d\rho(0)=d\phi(0),$  we have  for $\vert z\vert \leq 2\epsilon,$
\begin{eqnarray*}
\left\vert (d\rho_{\epsilon}-d\phi)(z)\right\vert&\leq&
\psi_{\epsilon}(z)\left\vert (d\rho-d\phi)(z)\right\vert+
\left\vert d\psi_{\epsilon}(z)\right\vert\cdot\left\vert
(\rho-\phi)(z)\right\vert\\
&\leq& \frac{1}{8}+\frac{C^{'}\epsilon^2 }{\epsilon},
\end{eqnarray*}
where $C^{'}$ is a finite constant.
Therefore there exists  $\epsilon_1>0$ such that for all
$0<\epsilon<\epsilon_1,$
\begin{equation}
\left\vert (d\rho_{\epsilon}-d\phi)(z)\right\vert\leq\frac{1}{4},\qquad
\vert z\vert \leq 2\epsilon.
\end{equation}
For any $0<\epsilon<\min\{\epsilon_0,\epsilon_1\}$  define
\begin{equation}
V^{'}:=\left\lbrace z\in\C^n:\ \rho_{\epsilon}(z)<0  \right\rbrace
\end{equation}
and let $V$ be the connected component of $V^{'}$ satisfying $P\in\partial
V.$
This, combined with (2.7) implies that $\vert d\rho_{\epsilon}(z)\vert
>\frac{1}{2}$ for $\vert z\vert \leq 2\epsilon.$ Since $0\leq \psi\leq 1$
and
$\rho_{\epsilon}(z)=\phi(z)$ for $\vert z\vert \geq 2\epsilon,$
we deduce from (2.6) that  $V$ is a $\mathcal{C}^2$ smooth subdomain of
$D.$

Now let
\begin{equation*}
U_3:=B(0,\epsilon),\qquad U_2:=U_3\cap T^{\R}_P,\qquad U_1:=U_3\cap T^{\C}_P
\end{equation*}
Then in virtue of (2.4)--(2.6), we see that Parts 1) and 2)  are satisfied
when  $\epsilon$ in (2.8) is sufficiently small.

We next turn to Part 3). Fix any $Q\in U_1$ and $z\in\partial V_Q,$ then there
are two cases. If $\vert z\vert \leq 2\epsilon,$
then by (2.7)
\begin{equation*}
\left\vert d_{z_1}\rho_{\epsilon} \right\vert
\geq \left\vert d_{z_1}\phi \right\vert
-\left\vert d \rho_{\epsilon}-d\phi \right\vert
>\frac{\vert z_1 +r\vert}{r}-\frac{1}{4}>\frac{1}{4}.
\end{equation*}
If  $\vert z\vert \geq 2\epsilon,$ then by (2.6)
\begin{equation*}
\left\vert d_{z_1}\rho_{\epsilon} \right\vert
 = \left\vert d_{z_1}\phi \right\vert
=\frac{\vert z_1 +r\vert}{r}>0.
\end{equation*}
Thus for any $Q=(0,z^{'})\in U_1$ the region  $V\cap \left\lbrace  (t,z^{'}),\  t\in \C
 \right\rbrace$   is a $\mathcal{C}^2$ smooth planar region contained in
 $D.$ Since for sufficiently small $\epsilon>0,$  $\partial D\cap U_3$ is
 a graph over $T^{\R}_P,$ a geometric argument shows that
 $V_Q$  is also a $\mathcal{C}^2$ smooth planar  simply connected region contained in
 $D.$ We see that one may assume that $V_Q$ is a domain.

 To complete Part 3) we still need to check that the family
 $(V_Q)_{Q\in U_1}$ satisfies 1) and 2) in Corollary 2.4.  Indeed,
 let $\rho_Q$ be the restriction  of $\rho$ on the complex line
 containing $V_Q.$ Clearly, the Hessian $d^2\rho_Q$ depends continuously on the
 parameter $Q\in U_1.$ This, combined with the proof of the geometric fact
 (see \cite[p. 325]{kr}) implies the remaining assertion of Part 3).

It  remains to establish Part 4). Also by  \cite[p. 325]{kr}, when
$\epsilon>0$ in (2.8) is sufficiently small, for any $z\in U_3\cap D$
there are a unique  point $\theta\in\partial D$ and
a unique point $\eta\in \partial V_Q$ such that
\begin{equation*}
\vert z-\theta\vert=\dist(z,\partial D)\qquad\text{and}\qquad
\vert z-\eta\vert=\dist(z,\partial V_Q).
\end{equation*}
Let $n_{\theta}$ (resp. $n_{\eta}$)  be the inward unit normal vector to
$\partial D$ (resp. $\partial V_Q$) at $\theta$ (resp. $\eta$).
Then  a geometric argument shows that the orthogonal projection
of the real line containing $n_{\eta}$ onto $V_Q$ passes through $z.$
Since $Q$ is close to $P,$ the angle between two vectors $z-\eta$ and
$n_{\eta}$ and the angle between $n_{\eta}$ and $n_{\theta}$ are
sufficiently small when $\epsilon$ in (2.8) is sufficiently small.
Thus the angle between two vectors $z-\eta$ and
$ z-\theta$  is
sufficiently small. Since $\vert d\rho(0)\vert=1$ and $\rho\in\mathcal{C}^2(U),$ it follows that
\begin{equation*}
\vert z-\theta\vert \leq \vert z-\eta\vert \leq C\vert z-\theta\vert
\end{equation*}
for some finite constant $C,$ which proves that
$\dist(z,\partial V_Q)\leq C\cdot \dist(z,\partial D).$

The second estimate of Part 4) can be proved in exactly the same way.
This completes the proof.
\end{proof}
\section{Estimates for the plurisubharmonic measures}
In this section we apply the result of the previous one  to
establish some inequalities concerning the
plurisubharmonic measures.
These estimates will be crucial for the proof of the main theorem.

\begin{prop}
Let $D$ be a bounded planar domain with $\mathcal{C}^2$ smooth boundary.
Then there is a constant $C$ with the following property:
For any   union $A$ of a finite number of  open connected arcs on $\partial D$,
one has
\begin{equation*}
\omega(z,A,D)\leq C\frac{\dist(z,\partial D)}{\dist(z,\partial D\setminus A)^2},\qquad z\in D\cup A.
\end{equation*}
\end{prop}
\begin{proof}
By Theorem 2.3 we know that there is a finite constant $C^{'}$ such that
\begin{equation*}
P_{D}(x,y)\leq C^{'}\frac{\dist (x,\partial D)}
{\vert x-y\vert^2},\qquad x\in D,\ y\in \partial D.
\end{equation*}
This, combined with identity (2.1), implies that
\begin{equation*}
\omega(z,A,D)\leq C^{'}\cdot\int\limits_{\partial D\setminus A}
\frac{\dist(z,\partial D) }{\vert z-\zeta\vert^2}d\sigma(\zeta),
\end{equation*}
where $d\sigma$ is the Lebesgue measure on $\partial D.$
We easily see  that the right side of the latter estimate is dominated by
\begin{equation*}
C\frac{\dist(z,\partial D)}{\dist(z,\partial D\setminus  A)^2}.
\end{equation*}
Hence the proof is  complete.
\end{proof}

   Observe that as in Theorem 2.3 and Corollary
 2.4, the constant $C$ in Proposition 3.1 depends only on $\sigma(\partial D\setminus A),$
$\diam(D)$ and the radius $r(D).$

\begin{prop}
Let $D\subset \C^n$  be a bounded  open set and let $A$ be an open set of $\partial
D$ such that $D$ is locally  $\mathcal{C}^2$ smooth on $A.$ Then for
any   set $K\Subset A,$  there is a finite constant $C=C_K$ such that
\begin{equation*}
\omega(z,A,D)  \leq C\cdot\dist(z,K),\quad
z\in D\cup A.
\end{equation*}
In particular, $\omega(\cdot,A,D)=0$  on $ A.$
\end{prop}
\begin{proof}
Since $\omega(\cdot,A,D)\leq \omega(\cdot,B,G)$ if $B\subset A$ and $G\subset D$ and
by using a compactness argument and applying Proposition 2.5,  we may
suppose without loss of generality that  $K\Subset A$ is the intersection of a  sufficiently small
ball centered  $U$ at $P$ with $A$ and $D$ is a  $\mathcal{C}^2$ smooth domain such that Proposition 2.5
is applicable in this context.  Namely, keeping    the notation in (2.2) and
(2.3), we assume without loss of generality that
$P= 0\in\C^n$ and  $\left(V_{Q}\right)_{Q\in U_1}$ is
 a   family of   $\mathcal{C}^2$ smooth planar simply connected
 domains satisfying 1) and 2) of Corollary 2.4.

 Observe that it suffices to prove the proposition for the case where $z$
 is sufficiently close  to $K.$ Now let $Q:=\pi^{\C}(z)$ and note that $z\in D_Q.$
 Then under this assumption, Part 4) of Proposition 2.5 gives a finite constant $C^{'}$
 such that
 \begin{equation}
\dist(z,K)\leq \dist(z,\partial V_Q)\leq C^{'} \cdot\dist(z,K).
 \end{equation}
Combining  Propositions 2.5, 3.1 and the remark made at the end of the proof of Proposition 3.1,
 we see that there is a finite constant $C^{''}$ such that
 \begin{equation}
\omega(z, A\cap \partial V_Q,V_Q)\leq C^{''}\cdot\dist (z,K\cap\partial V_Q).
 \end{equation}
Next observe that for any $u\in \mathcal{PSH}(D)$ with $u\leq 1$ on $D$ and $\hat{u}\leq 0$ on $A,$
 the following estimate
holds by the very definition
\begin{equation*}
u(z)\leq \omega(z, A\cap \partial V_Q,V_Q)  .
\end{equation*}
This, combined with (3.1) and (3.2)  implies that
\begin{equation*}
\omega(z,A,D)\leq  C^{'}C^{''}\cdot\dist (z,K ),
\end{equation*}
which completes the proof of the first desired estimate. The desired
identity $\omega(\cdot,A,D)=0$ on $A$ follows immediately
from this estimate. Hence the proof is finished.
\end{proof}
The next result  tells us that the definition of the
plurisubharmonic measure formulated  in Case II in Subsection 1.1 is well-defined.
\begin{prop}
Let $D\subset \C^n$  be an open set and let $A$ be an open subset of $\partial
D$ such that $D$ is locally   $\mathcal{C}^2$ smooth on $A.$
Then there is a function plurisubharmonic in $D$ which we denote by
$\omega(\cdot,A,D),$ with the following property:

Let $\left ( D_k \right )^{\infty}_{k=1}$ be a sequence
of relatively compact open sets $D_k \Subset
\C^n$  and $\left ( A_k\right )^{\infty}_{k=1}$  a sequence of open subsets of $A$  such that
\begin{itemize}
\item[(i)] $D_k\subset D_{k+1}$ and $\cup_{k=1}^{\infty} D_k=D;$
\item[(ii)]  $A_k\subset A_{k+1}$ and $A_k\subset \partial D\cap\partial D_k$
and $\cup_{k=1}^{\infty} A_k=A;$
 \item[(iii)]  for any point $\zeta\in A$ there is an open neighborhood
 $V=V_{\zeta}$ of $\zeta$ in $\C^n$ such that $V\cap D$ =$V\cap D_k$ for
 some $k.$
  \end{itemize}
Then
\begin{equation*} \omega(\cdot, A,  D)=
\lim\limits_{k\to\infty}
\omega(\cdot, A_k,  D_k)
  \qquad \text{on}\ D.
\end{equation*}
\end{prop}
\begin{proof} First observe that such  sequences $\left(D_{k}\right)_{k=1}^{\infty}$
and $\left(A_{k}\right)_{k=1}^{\infty}$ always exist. For example, one may
take $D_k:= D\cap B(0,k)$ and $A_k:= A\cap B(0,k),$ $k\in\N.$
Let  $\left(D^{'}_{k}\right)_{k=1}^{\infty}$  and
$\left(A^{'}_{k}\right)_{k=1}^{\infty}$
be  another  sequences
which verify (i)--(iii).
It is easy to see that the following limits of decreasing sequences
\begin{equation*}
u:=\lim\limits_{k\to\infty}
\omega(\cdot,  A_k,  D_k)\quad\text{and}\quad u^{'}:=\lim\limits_{k\to\infty}
\omega(\cdot,  A^{'}_k,  D^{'}_k)
\end{equation*}
exist and define two plurisubharmonic  functions in $D.$

Fix an $k$ and let $\zeta$ be any point in $ A_k$ and $K$ be any compact
  neighborhood of $\zeta$ in  $A_k.$ In virtue of (i)--(iii), there
are
a sufficiently large integer  $N$   and a bounded open   neighborhood $U$ of $K$ in $\C^n$
  such that $ U\cap A\subset A^{'}_n $ and $U\cap D^{'}_n=U\cap D$
   for any $n\geq N.$

Therefore,  applying Proposition 3.2 to the open set $D\cap U,$ we may find a  finite
constant $C$ such that
\begin{equation*}
 \omega(z,A^{'}_n, D^{'}_n )\leq  \omega(z,U\cap A, U\cap D )
\leq C\cdot \dist(z,K), \qquad z\in U\cap D,\ n\geq N.
\end{equation*}
 This implies that
\begin{equation*}
\widehat{u^{'}}(\zeta):=\limsup_{w\in D,\ w\to \zeta } u^{'}(w)=0.
\end{equation*}
Thus $\widehat{u^{'}}=0$ on $A_k$ and therefore
$ \omega(\cdot, A_k,  D_k)    \geq u^{'}.$ This implies
  that $u\geq u^{'}.$ Similarly, one gets  $ u^{'}\geq u.$
  Hence $u=u^{'}$ and the proof is finished.
\end{proof}

 One should mention that in virtue of
Proposition 3.3, Proposition 3.2 still holds when $D$ is
an arbitrary (not necessarily bounded) open set. An immediate consequence
of Proposition 3.3 is the following result.
\begin{prop}
Let $D\subset \C^n$  be an open set and let $A$ be an open set of $\partial
D$ such that $D$ is locally  $\mathcal{C}^2$ smooth on $A.$
Let $\left(A_{k}\right)_{k=1}^{\infty}$ be a sequence of open subsets
 of $\partial D$  such that $A_k\nearrow A$ as $k\nearrow \infty.$  Then
\begin{equation*}
\lim\limits_{k\to\infty}
\omega(\cdot, A_k,  D)
= \omega(\cdot, A,  D)\qquad \text{on}\ D.
\end{equation*}
\end{prop}
\begin{proof}
 It suffices to choose the sequence  $\left( D_k\right)^{\infty}_{k=1}$ with
 $D_k=D.$ Then
 the desired conclusion follows from Proposition 3.3.
\end{proof}
\begin{prop}
Let $D\subset \C^n$  be an open set and let $A$ be an open set of $\partial
D$ such that $D$ is locally   $\mathcal{C}^2$ smooth on $A.$ Then
for any $\delta>0,$ there is an open subset $T_{\delta}$ of $D$ such that
\\
1) $T_{\delta_1}\subset T_{\delta_2}$ for $0<\delta_1<\delta_2,$\\
2)
  $T_{\delta}\cup A$ is an open neighborhood of $A$ in $A\cup D,$
  \\
  3)
$\omega(z,A,D)-\delta
\leq \omega(z, T_{\delta},D) \leq\omega(z,A,D)$ for $ z\in
D,$  and\\
4) $\sup\limits_{T_{\delta}}\dist (\cdot,A)<\delta.$
\end{prop}
\begin{proof}
Fix a sequence $\left( A_k\right)^{\infty}_{k=1}$ of
open subsets of $A$ such that
\begin{itemize}
\item[(i)] $A_k\Subset A_{k+1}\subset A,$
\item[(ii)] $A_k\nearrow A$ as $k\nearrow\infty,$
\item [(iii)] $A_k$ consists of finite open connected components.
\end{itemize}
By Propositions 3.2 and  3.3, for any $k$ there is a finite constant $C_k>1$ such that
\begin{equation}
 \omega(z, A,D)\leq C_k \dist(z,A_k) .
\end{equation}
For an $\delta>0$ consider the following open subset of $D$
\begin{equation}
T_{\delta}:=\left\lbrace  z\in D:\   C_k \dist(z,A_k)<\delta\ \ \text{for some}\ \ k\in \N  \right\rbrace
\end{equation}
In virtue of (3.3)--(3.4) and (i)--(iii),   Part 1) and 2) are  proved.
Moreover
\begin{equation*}
\omega(z,A,D)\leq\delta \quad\text{and}\quad \dist(z,A)<\delta, \qquad z\in T_{\delta}
   ,
\end{equation*}
which  implies that
$    \omega(\cdot,A,D)-\delta  \leq \omega(\cdot, T_{\delta},D)             $ on $D.$

On the other hand,  by Part 2) and the definition of plurisubharmonic measures, we deduce that
  $\omega(\cdot, T_{\delta},D) \leq\omega(\cdot,A,D)$ on $D.$
Hence the proof is complete.
\end{proof}

The rest of  this section is devoted to   some applications of the previous results.

\begin{prop}
Let $D\subset \C^n, G\subset \C^m$  be two  domains and let $A$ (resp. $B$) be an
open set of $\partial
D$ (resp. $\partial G$) such that $D$ (resp. $G$) is locally   $\mathcal{C}^2$
smooth on $A$ (resp. $B$).
Put $X:=\X(A,B;D,G)$ and $ \widehat{X}^{\text{o}}:=
     \widehat{\X}^{\text{o}}(A,B;D,G).$ Then\\
\noindent 1) for any finite subset
     $M\subset \widehat{X}^{\text{o}},$ there are open sets
     $T\subset D,$   $S\subset G$ and
  $0<\epsilon<1$  such that
\begin{equation*}
M\subset \left\lbrace (z,w)\in D\times G:\ \omega(z,T,D) +\omega(w, S,G) <1-\epsilon
\right\rbrace
\subset  \widehat{X}^{\text{o}};
\end{equation*}
\noindent 2) the open set $ \widehat{X}^{\text{o}}$ is connected;

\noindent 3) $X\subset \widehat{X}.$
\end{prop}
\begin{proof}
Fix an $\epsilon>0$ such that
\begin{equation*}
\omega(z,A,D) +\omega(w, B,G) <1-2\epsilon,\qquad (z,w)\in M.
\end{equation*}
Applying Proposition 3.5, we find two open sets  $T\subset D,$   $S\subset G$
of the form (3.4) such that
\begin{eqnarray*}
\left\vert \omega(z,A,D) -\omega(z, T,D)\right\vert &<&\frac{\epsilon}{2},\qquad z\in D,  \\
\left\vert \omega(w,B,G) -\omega(w, S,G)\right\vert &<&\frac{\epsilon}{2},\qquad w\in G  .
\end{eqnarray*}
 Therefore,
\begin{equation*}
M\subset   \left\lbrace (z,w)\in D\times G:\ \omega(z,T,D) +\omega(w, S,G) <1-\epsilon
\right\rbrace\subset  \widehat{X}^{\text{o}},
\end{equation*}
 which finishes Part 1).

To prove Part 2) let $(z_1,w_1)$ and $(z_2,w_2)$ be two arbitrary points in $ \widehat{X}^{\text{o}}.$
Put $M:=\left\lbrace (z_1,w_1), (z_2,w_2)\right\rbrace .$ By
Part 1) there are open sets
     $T\subset D,$   $S\subset G$ and
  $0<\epsilon<1$  such that
\begin{equation*}
M\subset \left\lbrace (z,w)\in D\times G:\ \omega(z,T,D) +\omega(w, S,G) <1-\epsilon
\right\rbrace
\subset  \widehat{X}^{\text{o}}
.
\end{equation*}
Since the set $\left\lbrace (z,w)\in D\times G:\ \omega(z,T,D) +\omega(w, S,G) <1-\epsilon
\right\rbrace$ is connected (see, for example, Lemma 4 in \cite{jp2}), the desired
conclusion of Part 2) follows.

Part 3) holds by applying Proposition 3.2 and taking into account  the remark made just
after Proposition 3.3.  Hence the proof is complete.
\end{proof}

The next result tells us that  the  open set $
\widehat{X}^{\text{o}}$
is still connected in the following mixed situation.
\begin{prop}
 Let $D\subset \C^n,$ $G\subset \C^m$ be two domains and  let
  $A\subset D$ and   $B$ is an open subset of $\partial G.$
  Assume that $A$ is locally pluriregular  and $G$ is locally  $\mathcal{C}^2$ smooth
   on $B.$
    Let $X:= \X(A,B;D,G)$   and      $ \widehat{X}^{\text{o}}:=
     \widehat{\X}^{\text{o}}(A,B;D,G).$
Then\\
1)  for any finite subset
     $M\subset \widehat{X}^{\text{o}},$ there are an open set
     $S\subset G$ and a number
  $0<\epsilon<1$  such that
\begin{equation*}
M\subset \left\lbrace (z,w)\in D\times G:\ \omega(z,A,D) +\omega(w, S,G) <1-\epsilon
\right\rbrace
\subset  \widehat{X}^{\text{o}}
;
\end{equation*}
 \noindent 2)    the open set $ \widehat{X}^{\text{o}}$ is connected;\\
 3) $X\subset \widehat{X}.$
\end{prop}
\begin{proof}  We proceed as in
the proof of Proposition 3.6
making the obviously necessary changes.
Hence, the proof is complete.
\end{proof}

The last result of this section studies the level sets of
plurisubharmonic measures.

\begin{prop}
Let $D\subset \C^n$  be an open set and let $A$ be an open set of $\partial
D$ such that $D$ is locally  $\mathcal{C}^2$ smooth on $A.$
 For any $0<\epsilon <1$ let
\begin{equation*}
D_{\epsilon}:=\left\lbrace z\in D:\ \omega(z,A,D)<1-\epsilon\right\rbrace.
\end{equation*}
 Then\\
 1)
\begin{equation*}
\lim\limits_{\epsilon\to 0}
\omega(\cdot, A,  D_{\epsilon})
= \omega(\cdot, A,  D)\quad \text{on}\ D\ \text{and}\
\omega(\cdot,A,D_{\epsilon})=\frac{\omega(\cdot,A,D)}{1-\epsilon}\quad \text{on} \
D_{\epsilon};
\end{equation*}
2) if $z\in D$ verifies $\omega(z,A,D)<1,$ then,   for any $0<\epsilon<1-\omega(z,A,D),$
the connected component of $ D_{\epsilon}$ which contains $z$ is
locally  $\mathcal{C}^2$ smooth on a nonempty  open subset of $A;$
\\
3) for   any $0<\epsilon_0<1,$ there is an open
neighborhood $U$ of $A$ in $D\cup A$
such that for all $\epsilon \leq 1-\epsilon_0$ there exists exactly one
connected component of $D_{\epsilon}$ containing $U\cap D$ and
$\omega(\cdot,A,D_{\epsilon})<\epsilon_0$ on $U.$
\end{prop}
\begin{proof}
For  $k\in \N$ let $A_{k}:=A.$ It suffices to check that
the sequences  $\left( D_{\frac{1}{k}}\right)^{\infty}_{k=N}$
and  $\left( A_k\right)^{\infty}_{k=N}$ satisfy the properties (i)-(iii)
of Proposition 3.3 for a sufficiently large positive integer $N.$
 Observe that the only nontrivial verification is
for (iii).
But this follows immediately from an application of Proposition 3.2.
Hence the first identity of Part 1) is proved.
To verify the second one, observe by the very definition that
$\omega(\cdot,A,D_{\epsilon})\geq\frac{\omega(\cdot,A,D)}{1-\epsilon}$ on
$D_{\epsilon}.$ On the other hand, consider the function
\begin{equation*}
 u(z):=
\begin{cases}
\text{max}\left\lbrace \omega(z,A,D),(1-\epsilon)\omega(z,A,D_{\epsilon})  \right\rbrace,
  & z\in    D_{\epsilon},\\
    \omega(z,A,D)   , & z \in D\setminus D_{\epsilon}.
\end{cases}
\end{equation*}
Clearly, $u\in\mathcal{PSH}(D)$ and $u\leq 1$ on $D.$ By Proposition 3.2,
$\hat{u}=0$ on $A.$ Thus $u\leq \omega(\cdot,A,D)$ which completes
the proof of the second identity of Part 1).

 By Part 1),
$\omega(z,A,D_{\epsilon})<1$ for all $0<\epsilon<1-\omega(z,A,D).$
This proves Part 2).

For Part 3) it suffices  to show that for every point $\zeta\in A$ there
is a neighborhood $U$ of $\zeta$ in $D\cup A$ which possesses the required
properties.
By Proposition 3.2 one may find an open
neighborhood $\mathcal{U}$ of $\zeta$ such that for some finite constant $C_1,$
\begin{equation*}
\omega(\cdot, A,  D)
\leq C_1\cdot\dist(z,\mathcal{U}\cap A)<\epsilon_0\qquad \text{on}\ \mathcal{U}\cap D.
\end{equation*}
This shows that $\mathcal{U}\cap D\subset D_{1-\epsilon_0}.$
We now choose a  relatively compact neighborhood $U$ of $\zeta$ such that
$U\Subset \mathcal{U}.$
Then  applying Proposition 3.2 and  shrinking $U,$ if necessary, we also have
\begin{equation*}
\omega(\cdot, A,  D_{1-\epsilon_0})
\leq C_2\cdot\dist(z,A)<\epsilon_0\qquad \text{on}\ U\cap
D,
\end{equation*}
which completes the last part of the proposition.
\end{proof}
\section{A mixed  cross theorem and  two quantitative  cross theorems}
The main result of this section  is the following  mixed cross theorem.
\begin{thm}
 Let $D\subset \C^n$ be a bounded pseudoconvex domain, $G\subset \C^m$  a
 domain,
  $A\subset D,$ and   $B\subset\partial G.$
  Assume that $A$ is a locally pluriregular relatively compact subset of $D$ and
  $A=\bigcup\limits_{k=1}^{\infty} A_k$ with $A_k$ locally pluriregular  compact subsets of $D$
  and that
   $B$ is an open subset of  $\partial G$ such that $G$ is locally  $\mathcal{C}^2$ smooth
   on $B.$
    Let $X:= \X(A,B;D,G)$, $X^{\text{o}}:= \X^{\text{o}}(A,B;D,G),$
    $\widehat{X}:= \widehat{\X}(A,B;D,G),$ and      $ \widehat{X}^{\text{o}}:=
     \widehat{\X}^{\text{o}}(A,B;D,G).$

Let $\mathcal{C}_s(X)$ be the space  of all functions defined on $X$ such
that
\begin{itemize}
\item[(i)]  $f$ is locally bounded on $X;$
\item[(ii)] for any $z\in A,$  $f(z,\cdot)\in
\mathcal{C}(G\cup B).$
\end{itemize}
      1) Then for  any   function $f\in \mathcal{C}_s(X)\cap
      \mathcal{O}_s(X^{\text{o}})$
        there is a unique function
$\hat{f}\in \mathcal{C}(\widehat{X})\cap\mathcal{O}(\widehat{X}^{\text{o}})$
such that $\hat{f}=f$ on $X .$\\
 2) If, moreover, there is a set $B^{'}\subset \partial G$ such that
\begin{itemize}
\item[$(\text{i}^{'})$] $f$ is locally bounded on $A\times (G\cup B^{'}),$
\item[$(\text{ii}^{'})$] for any $z\in A,$  $f(z,\cdot)\in
\mathcal{C}(G\cup B^{'}),$
\item[$(\text{iii}^{'})$]  $\omega(\cdot,B,G)\in\mathcal{C}(G\cup B^{'}),$
\end{itemize}
then $\hat{f}$ extends continuously to every point $(z,\eta)\in
\overline{\widehat{X}}\cap(D\times B^{'}).$
\end{thm}

A remark is in order. Under the hypothesis of Theorem 4.1, it follows from
Part 3) of Proposition 3.7  that $X\subset\widehat{X}.$
\begin{proof} First we prove Part 1).
We argue as in the proof of Theorem 3.5.1  in \cite{jp1}. For the sake of
completeness, we give here a sketchy proof.
Fix an $f\in\mathcal{C}_s(X)\cap\mathcal{O}_s(X^{\text{o}}).$

\smallskip

\noindent{\bf Step I:} {\it Reduction to the case where $D$ is strongly
pseudoconvex, $A$ is a locally pluriregular  compact subset of $D$
and $\vert f\vert$ is bounded on $X.$}

\smallskip

One proceeds as in the first and second step in that
proof. More precisely, let  $(G_k)_{k=1}^{\infty}$ be an exhausting sequence of $G$
which verifies the properties (i)--(iii) in Proposition 3.3 (with $B$ instead of $A$).
Let $B_k:=B\cap\partial G_k.$ Since $D$ is a domain of holomorphy, we may
find an exhausting sequence  $(D_k)_{k=1}^{\infty}$ of relatively compact,
strongly pseudoconvex subdomains of $D$ with $A_k\subset D_k\nearrow D.$

By reduction assumption, for each $k$ there exists an
$\hat{f}_k\in\mathcal{C}\left(\widehat{\X}(A_k,B_k; D_k,G_k)\right)\cap
\mathcal{O}\left(\widehat{\X}^{\text{o}}(A_k,B_k; D_k,G_k)\right)$
such that $\hat{f}_k=f$ on $\X(A_k,B_k; D_k,G_k).$

By Theorem 2.1  and  Proposition 3.7 and taking into account that
$f_{k+1}=f_k=f$ on   $\X(A_k,B_k; D_k,G_k)$, one can
show that $f_{k+1}=f_k=f$ on   $\widehat{\X}(A_k,B_k; D_k,G_k).$
On the other hand, by Proposition 3.3, one gets $\widehat{\X}(A_k,B_k; D_k,G_k) \nearrow \widehat{\X} $
as $k\nearrow \infty.$
Therefore, we may glue $f_k$ together to obtain
 a function
$\hat{f}\in \mathcal{C}(\widehat{X})\cap\mathcal{O}(\widehat{X}^{\text{o}})$
such that $\hat{f}=\hat{f}_k=f$ on $\X(A_k,B_k;D_k,G_k) .$ Thus $\hat{f}=f$ on $X.$
The uniqueness of such an extension $\hat{f}$  follows from
Theorem 2.1 and Proposition 3.7.
This completes Step I.

\smallskip

\noindent{\bf Step II:} {\it The case where
$D$ is strongly
pseudoconvex, $A$ is a locally pluriregular  compact subset of $D$
and  $\vert f\vert \leq 1$   on $X.$}

The key observation is that we are still able to apply the classical
method of doubly orthogonal bases of Bergman type
(see for example \cite{ns}, \cite{nz} for a systematic study of this method).

Next one observes that Lemma 3.5.10 in \cite{jp1} is still valid in the present context.
Look at  Step 3 in that proof. In the sequel, we will use the notations from \cite{jp1}.

Let $\mu:=\mu_{A,D},$ $H_0:= L^2_{h}(D),$ $H_1:=$ the closure of
$H_0|_{A}$ in $L^2(A,\mu)  $ and let $(b_k)_{k=1}^{\infty}$ be the basis
from Lemma 3.5.10  in \cite{jp1}, $\nu_k:=\Vert b_k\Vert_{H_0},\ k\in\N,$
and $\nu_k\nearrow \infty.$

For any $w\in B,$ we have $f(\cdot,w)\in H_0$ and   $f(\cdot,w)|_{A}\in
H_1.$ Hence
\begin{equation}
f(\cdot,w)=\sum\limits_{k=1}^{\infty} c_k(w)b_k,
\end{equation}
where
\begin{equation}
c_k(w)=\frac{1}{\nu_k^2}\int\limits_{D} f(z,w)\overline{b_k(z)}
d\Lambda_{2n}(z)=\int\limits_{A} f(z,w)\overline{b_k(z)}
d\mu(z),\qquad k\in\N.
\end{equation}
Taking  the hypothesis $\vert f\vert \leq 1$  on $X$ and
$f\in\mathcal{C}_s(X)\cap\mathcal{O}_s(X)$
into account and applying   the Lebesgue's Dominated Convergence
Theorem, we see that the formula
\begin{equation}
\widehat{c_k}(w):=\int\limits_{A} f(z,w)\overline{b_k(z)}
d\mu(z),\qquad w\in G\cup B,\  k\in\N;
\end{equation}
defines a bounded  function which is holomorphic in $G.$
Moreover, by (4.2)--(4.3)   it follows that
\begin{equation}
\lim_{w\in G,\ w\to \eta }\widehat{c_k}(w)  =\widehat{c_k}(\eta)=c_k(\eta),\qquad \eta\in B .
\end{equation}
Thus $ \widehat{c_k}\in \mathcal{C}(G\cup B)\cap   \mathcal{O}(G).$

 Observe that as in
\cite{jp1} and using (4.2)--(4.4), we  obtain the following estimates
\begin{eqnarray*}
\frac{\log{\vert \widehat{c_k}\vert} }{\log{\nu_k}}&\leq & \frac{\sqrt{\mu(A)}
}{\log{\nu_k}},\qquad k\in\N,\\
 \limsup_{w\in G,\ w\to \eta }\frac{\log{\vert \widehat{c_k}(\eta)\vert} }{\log{\nu_k}}&\leq &
\frac{\log{\sqrt{\Lambda_{2n}(D)}    } }{\log{\nu_k}}-1, \qquad \eta\in B,\
k\in\N.
\end{eqnarray*}
This shows that for any $\epsilon>0,$ there is a sufficiently large $N$ such that
for all $k\geq N,$
\begin{equation}
\frac{\log{\vert \widehat{c_k}\vert}
}{\log{\nu_k}}\leq\omega(\cdot,B,G)+\epsilon-1 \qquad\text{on}\ G.
\end{equation}

Take a compact $K\Subset D$ and let $\alpha>\max\limits_{K}
h^{\ast}_{A,D}$ and $\epsilon>0$ so small such that $\alpha+2\epsilon<1.$
Consider the open set
\begin{equation*}
G_K:=\left\lbrace w\in G:\ \omega(\cdot,B,G)<1-\alpha-2\epsilon
\right\rbrace.
\end{equation*}
By (4.5) there is a constant $C^{'}(K)$ such that
\begin{equation}
\Vert  \widehat{c_k}\Vert\leq C^{'}(K)\nu_k^{\omega(\cdot,B,G)+\epsilon-1}
\leq  C^{'}(K)\nu_k^{-\alpha-\epsilon},\qquad k\geq 1.
\end{equation}
Now we wish to show that
\begin{equation}
\sum\limits_{k=1}^{\infty}\widehat{c_k}(w)b_k(z)
\end{equation}
converges locally uniformly in $\widehat{X}^{\text{o}}.$
Indeed, by (4.6) and Lemma 3.5.10 in \cite{jp1},
\begin{equation*}
\sum\limits_{k=1}^{\infty}\Vert \widehat{c_k}\Vert_{G_K}\Vert
b_k\Vert_K\leq \sum\limits_{k=1}^{\infty}  C^{'}(K)\nu_k^{-\alpha-\epsilon}
C(K,\alpha)\nu_k^{\alpha }
\leq  C^{'}(K)C(K,\alpha)\sum\limits_{k=1}^{\infty}\nu_k^{-\epsilon
}<\infty,
\end{equation*}
which gives the normal convergence on $K\times G_K.$
Moreover,  by Proposition 3.2 one gets  $B\subset \partial G_K.$ Therefore,
the previous argument also shows that
the series in (4.7) converges normally on $K\times (G_K\cup B).$
Since the compact set $K\Subset D$ and $\epsilon>0$ are arbitrary,
the series in (4.7) converges uniformly on compact subsets of
$\widehat{X}.$  Let $\hat{f}$
be the sum limit. Then obviously $\hat{f}\in\mathcal{C}(\widehat{X})\cap\mathcal{O}(\widehat{X}^{\text{o}}).$
Taking (4.1), (4.4) and (4.7) into account, it follows that
 \begin{equation*}
   \hat{f}=f\qquad \text{on} \ D\times B.
\end{equation*}
 Consequently, an application of
 Theorem 2.1  gives that $\hat{f}=f$ on $ X.$
This completes the proof of Part 1).

We now turn to Part 2) using the proof of Part 1).
Observe that by (4.3) and $(\text{i}^{'})$--$(\text{ii}^{'}),$
$ \widehat{c_k}\in \mathcal{C}(G\cup B^{'})\cap   \mathcal{O}(G).$
Next fix an $\eta\in B^{'}$ and $z\in D.$ We use hypothesis $(\text{iii}^{'})$ in order to choose
 $\epsilon >0$
and a  compact  neighborhood $K$ of $z$ such that $K\times (G_K\cup\{\eta\})$ is a
neighborhood of $(z,\eta)$ in $\overline{\widehat{X}}\cap (D\times B^{'}).$
The rest of the proof follows essentially along the same lines as that of
Part 1). This completes the proof of Part 2).
\end{proof}

The last two results of this section give quantitative versions of the
classical cross theorem (cf. Theorem 1).
\begin{thm}
 Let $D\subset \C^n,$ $G\subset \C^m$
 be   bounded  domains and  let
  $A\subset D,$   $B\subset G$ be locally pluriregular sets.
  Assume that $D$ is pseudoconvex and  $A$ is  of the form
  $A=\bigcup\limits_{k=1}^{\infty} A_k$ with $A_k$  compact subset of $D.$
Let $X:= \X(A,B;D,G)$ and
    $\widehat{X}:= \widehat{\X}(A,B;D,G).$
    Then for any  $f\in  \mathcal{O}_s(X )$
        there is a unique function
$\hat{f}\in  \mathcal{O}(\widehat{X}  )$
such that $\hat{f}=f$ on $X .$ Moreover, if $\vert
f\vert_X<\infty$  then
\begin{equation}
\vert \hat{f}(z,w)\vert\leq \vert f\vert_{A\times B}^{1-\omega(z,A,D)-\omega(w,B,G)}
\vert f\vert_X^{\omega(z,A,D)+\omega(w,B,G)},\qquad
(z,w)\in\widehat{X}.
\end{equation}
\end{thm}
\begin{proof}
We proceed in two steps.

\smallskip

\noindent {\bf Step 1:}
{\it Proof of the equality $\vert \hat{f}\vert_{\widehat{X}}=\vert f\vert_X.$}

In order to reach a contradiction suppose  that there is a point
$(z_0,w_0)\in\widehat{X}$ such that $\vert\hat{f}(z_0,w_0)\vert>\vert
f\vert_X.$ Put $\alpha:=\hat{f}(z_0,w_0)$ and consider the function
\begin{equation}
g(z,w):=\frac{1}{f(z,w)-\alpha},\qquad  (z,w)\in X.
\end{equation}
Clearly, $g\in\mathcal{O}_s(X).$ Hence by Theorem 3.5.1 in \cite{jp1}, there is exactly
one function  $\hat{g}\in\mathcal{O}(\widehat{X})$ with $\hat{g}=g$ on
$X.$ Therefore, by (4.9) we have on $X:$ $g(f-\alpha)\equiv 1.$ Thus
$\hat{g}(\hat{f}-\alpha)\equiv 1$ on $\widehat{X}.$ In particular,
\begin{equation*}
0=\hat{g}(z_0,w_0)(\hat{f}(z_0,w_0)-\alpha)= 1,
\end{equation*}
which is a contradiction.
Hence the inequality $\vert \hat{f}\vert_{\widehat{X}}\leq\vert f\vert_X$ is
proved. The converse inequality is trivial since $X\subset\widehat{X}$ (see, for example, \cite{jp1}).  Thus Step 1 is complete.

\smallskip

\noindent {\bf Step 2:}
{\it Proof of   inequality (4.8).}

Fix now $(z_0,w_0)\in\widehat{X}.$
For every  $\eta\in B,$ we have
\begin{equation*}
\vert f(\zeta,\eta)\vert \leq \vert f\vert_{A\times
B},\ \zeta\in A \qquad\text{and}\qquad \vert f(z,\eta)\vert \leq \vert f\vert_{X},\ z\in
D.
\end{equation*}
Therefore, Two-Constant Theorem (Theorem 2.2) implies that
\begin{equation}
\vert f(z,\eta)\vert \leq \vert f\vert_{A\times
B}^{1-\omega(z,A,D)} \vert f\vert_{X}^{\omega(z,A,D)},\qquad z\in D,\
\eta\in B.
\end{equation}
Consider the function $\hat{f}(z_0,\cdot)\in\mathcal{O}(G_{z_0}),$
 where
 \begin{equation*}
 G_{z_0}:=\left\lbrace  w\in G: \omega(w,B,G)<1-\omega(z_0,A,D)\right\rbrace.
 \end{equation*}
 Observe that $\vert \hat{f}(z_0,\cdot)\vert_{G_{z_0}}\leq \vert f\vert_{X}$ and
$\omega(w,B,G_{z_0})=\frac{\omega(w,B,G)}{1-\omega(z_0,A,D)}.$ Consequently,
 using (4.10) and applying the Two-Constant Theorem  to the function $\hat{f}(z_0,\cdot),$
(4.8) for $(z_0,w_0)$ follows. Hence the theorem is proved.
\end{proof}

\begin{thm}
Let  $D\subset\C^{n},\ G\subset \C^m$  be   domains and
 let $A$ (resp. $B$) be an  open subset of the boundary $\partial D$ (resp.
  $\partial G$).
Suppose in addition  that $D$ (resp. $G$) is locally  $\mathcal{C}^2$ smooth on
  $A$ (resp. $B$) and $D$ is pseudoconvex.
   Put $X:=\X(A,B;D,G),$
$\widehat{X}:=\widehat{\X}(A,B;D,G)$  and $\widehat{X}^{\text{o}}:=\widehat{\X}^{\text{o}}(A,B;D,G).$
Then for any   function
$f\in\mathcal{C}(\widehat{X})\cap \mathcal{O}(\widehat{X}^{\text{o}})$
the following inequality holds
\begin{equation}
\vert f(z,w)\vert\leq \vert f\vert_{A\times B}^{1-\omega(z,A,D)-\omega(w,B,G)}
\vert f\vert_X^{\omega(z,A,D)+\omega(w,B,G)},\qquad
(z,w)\in\widehat{X}.
\end{equation}
\end{thm}
\begin{proof}
Fix a point $(z_0,w_0)\in\widehat{X},$ an arbitrary number $\epsilon>0$
and let $\delta>0.$

 By Proposition 3.5 one may find an open set
$T_{\delta}\subset D$ such that
\begin{equation}
\omega(z,A,D)-\delta\leq \omega(z,T_{\delta},D)\leq
\omega(z,A,D),\qquad z\in D.
\end{equation}
By Proposition 3.3 one may find a (not necessarily pseudoconvex) bounded subdomain  $G_{\delta}$ of $G$
such that $\overline{G_{\delta}}\Subset G\cup B,$
$G_{\delta}$ is locally $\mathcal{C}^2$ smooth on the open subset
$\partial G_{\delta}\cap B$ of $B$   and
\begin{equation}
0\leq \omega(w_0, \partial G_{\delta}\cap B      ,G_{\delta})- \omega(w_0,B,G)<\delta.
\end{equation}
Since $f\in\mathcal{C}(\widehat{X}),$ there is an open subset
$A_{\delta}$ of $T_{\delta}$ such that $A\cup A_{\delta}$ is
an open neighborhood of $A$ in $A\cup D$ and moreover
\begin{equation}
\vert f(z,w) \vert \leq \vert f\vert_X+\epsilon, \qquad z\in A_{\delta},\
w\in G_{\delta}.
\end{equation}
It is also clear from (4.12) and the above properties of $A_{\delta}$  that
\begin{equation}
\omega(z,A_{\delta},D)-\delta\leq \omega(z,T_{\delta},D)\leq
\omega(z,A_{\delta},D),\qquad z\in D.
\end{equation}

Let $D_{\delta}$ be a strongly pseudoconvex subdomain of $D$ such that
$D_{\delta}\Subset D$ and
\begin{equation}
0\leq \omega(z_0,  A_{\delta}\cap D_{\delta}      ,D_{\delta})- \omega(z_0,A_{\delta},D)<\delta.
\end{equation}

Since $G_{\delta}$ is locally $\mathcal{C}^2$ smooth on the open subset
$\partial G_{\delta}\cap B$ of $B$ and
$f\in\mathcal{C}(\widehat{X}),$ one may find an open subset
$B_{\delta}$ of $G_{\delta}$ such that $B\cup B_{\delta}$ is
an open neighborhood of $\partial G_{\delta}\cap B     $ in $(\partial G_{\delta}\cap B) \cup
G_{\delta} $ and moreover
\begin{equation}
\vert f(z,w) \vert \leq \vert f\vert_X+\epsilon, \qquad z\in D_{\delta},\
w\in B_{\delta}.
\end{equation}
By taking the intersection of $B_{\delta}$ with the level open set given by
Proposition 3.5 with respect to the open set $G_{\delta},$
one may assume that
\begin{equation}
 \omega(w,\partial G_{\delta}\cap B ,G_{\delta}) -\delta\leq\omega(w,B_{\delta},G_{\delta})\leq
\omega(w,\partial G_{\delta}\cap B ,G_{\delta}),\qquad w\in G_{\delta}.
\end{equation}
Consider  following crosses
\begin{equation*}
X_{\delta}:=\X(A_{\delta}\cap D_{\delta},B_{\delta};D_{\delta},G_{\delta})
\qquad\text{and}\qquad \widehat{X}_{\delta}:=
\widehat{\X}(A_{\delta}\cap D_{\delta},B_{\delta};D_{\delta},G_{\delta}).
\end{equation*}
By Theorem 1, there is a function $f_{\delta}\in\mathcal{O}(\widehat{X}_{\delta})$
such that $f_{\delta}=f$ on $X_{\delta}.$

If one chooses $\delta $ such that
$0<10\delta<1-\omega(z_0,A,D)-\omega(w_0,B,G),$
then it follows from (4.12), (4.13), (4.15), (4.16) and (4.18) that
\begin{equation*}
(z_0,w_0)\in\left\lbrace (z,w)\in D_{\delta}\times H_{\delta}:\quad
\omega(z,A_{\delta}\cap D_{\delta},D_{\delta})+\omega(w,B_{\delta},H_{\delta})<1-5\delta\right\rbrace
\subset \widehat{X}_{\delta},
\end{equation*}
where $H_{\delta}$ is the connected component of $G_{\delta}$ containing
$w_0.$

In addition we recall that $f_{\delta}=f$  on $X_{\delta}.$
Therefore, $f(z_0,w_0)=f_{\delta}(z_0,w_0).$ Consequently, applying
Theorem 4.2 and taking (4.14) and (4.17) into account,
we deduce that $\vert f(z_0,w_0)\vert \leq \vert f\vert_X+\epsilon.$
Since $\epsilon>0$ and $(z_0,w_0)\in\widehat{X}$ are arbitrary, it follows
that $\vert f\vert_{\widehat{X}}\leq \vert f\vert.$ The converse
inequality is trivial as $X\subset\widehat{X}$ by Part 3) of Proposition
3.6.
Thus we have shown that $\vert f\vert_{\widehat{X}}= \vert f\vert.$

Therefore, arguing as in Step 2 of Theorem 4.2 and applying
the second identity of Part 1) of Proposition 3.8, inequality (4.11)
follows.
\end{proof}
\section{Proof of the Main theorem for $N=2$}
In this section we simplify the notation and  rephrase the Main Theorem for the case $N=2$  as follows.
\begin{thm}
Let  $D\subset\C^{n},\ G\subset \C^m$  be  pseudoconvex domains and
 let $A$ (resp. $B$) be an  open subset of the boundary $\partial D$ (resp.
  $\partial G$).
Suppose in addition  that $D$ (resp. $G$) is locally  $\mathcal{C}^2$ smooth on
  $A$ (resp. $B$).
   Put $X:=\X(A,B;D,G),$   $X^{\text{o}}:=\X^{\text{o}}(A,B;D,G),$
$\widehat{X}:=\widehat{\X}(A,B;D,G)$  and $\widehat{X}^{\text{o}}:=\widehat{\X}^{\text{o}}(A,B;D,G).$
Then for any   function $f\in\mathcal{C}(X)\cap\mathcal{O}_s(X^{\text{o}}),$ there is a unique function
$\hat{f}\in\mathcal{C}(\widehat{X})\cap \mathcal{O}(\widehat{X}^{\text{o}})$
such that $\hat{f}=f$ on $X.$ Moreover,
\begin{equation}
\vert \hat{f}(z,w)\vert \leq \vert f\vert_{A\times B}^{1-\omega(z,A,D)-\omega(w,B,G)}
\vert f\vert_X^{\omega(z,A,D)+\omega(w,B,G)},\qquad
(z,w)\in\widehat{X}.
\end{equation}
\end{thm}
\begin{proof}
We proceed by several steps. First observe that by Theorem 2.1 and Part 3) of
Proposition 3.6, the function
$\hat{f}$ is uniquely determined (if exists).

\smallskip

\noindent{\bf Step 1:} {\it Reduction to the case where $D$ and $G$ are
bounded pseudoconvex domains}

\smallskip

 \noindent{\it Proof of Step 1.} Fix any sequences of bounded pseudoconvex
 subdomains $\left (D_k\right)_{k=1}^{\infty}$
 (resp. $\left (G_k\right)_{k=1}^{\infty}$)    of $D$ (resp. $G$)
  such that the sequences $\left (D_k\right)_{k=1}^{\infty}$ and $\left (A_k\right)_{k=1}^{\infty}$
(resp. $\left (G_k\right)_{k=1}^{\infty}$ and $\left
(B_k\right)_{k=1}^{\infty}$)
satisfy (i)--(iii) of Proposition 3.3, where  $A_k:= A\cap\partial
D_k$ and $G_k:=B\cap\partial G_k.$
Let
\begin{equation*}
X_k:=\X(A_k,B_k;D_k,G_k)\subset X
\end{equation*}
and note that $\widehat{X_k}\nearrow\widehat{X}$ by Propositions 3.3.

Let  $f\in\mathcal{C}(X)\cap\mathcal{O}_s(X^{\text{o}}) $ be given.
Clearly, $f\in\mathcal{C}(X_k).$ Therefore, by the
reduction assumption, for each $k$ there exists an $\hat{f}_k\in  \mathcal{C}(\widehat{X_k})\cap
\mathcal{O}(\widehat{X_k}^{\text{o}})$ with $\hat{f}_k=f$ on $X_k.$ By Theorem 2.1 and Proposition 3.6,
$ \hat{f}_{k+1}=\hat{f}_k    $ on $\widehat{X_k}.$ Therefore, gluing the
$\hat{f}_k$'s, we obtain an $\hat{f}\in \mathcal{C}(\widehat{X})\cap  \mathcal{O}(\widehat{X}^{\text{o}})   $
 with $\hat{f}=f$ on $X.$
To reduce  estimate (5.1) to the case where $D$ and $G$ are bounded pseudoconvex domains,
 we proceed in the same way as above.
 This completes Step 1.
\hfill $\square$

\smallskip

{\it From now on we assume that the hypothesis of Step 1 is fulfilled. }

We introduce a new terminology.  A subset $\mathcal{A}$ of an open subset $A$ of $\partial D$
 is said to be  {\it a ball in $A$ with center $\zeta$ and radius $r$} if
$\mathcal{A}=B(\zeta,r)\cap A$  for a point $\zeta\in A$ and a positive  number $r$
verifying  $2r<\dist(\zeta,\partial A).$
 Moreover, for a ball $\mathcal{A}$ in $A$ and a number $0<\lambda\leq 2,$ $\lambda\mathcal{A}$
 denotes the open set
$ B(\zeta,\lambda r)\cap \partial D.$

\smallskip

\noindent{\bf Step 2:} {\it We keep the hypothesis of Theorem 5.1 and assume in addition that
 $G$ is a Jordan
planar simply connected domain. Then we shall prove the following local version of Theorem 5.1:}

\smallskip

{\it For any point $P\in A,$ there is  a ball
$\mathcal{A}$ in $A$ with center  $P$ such that
 the following property holds:  For any   function
 $f\in\mathcal{C}(X)\cap\mathcal{O}_s(X^{\text{o}}),$ there
  is a unique function
$\hat{f}\in\mathcal{C}(\widehat{X}_{\mathcal{A}})\cap \mathcal{O}(\widehat{X}^{\text{o}}_{\mathcal{A}})$
such that $\hat{f}=f$ on $X_{\mathcal{A}},$
where
 \begin{equation*}
 X_{\mathcal{A}}:=\X(\mathcal{A},B; D,G), \
 \widehat{X}_{\mathcal{A}}:=\widehat{\X}(\mathcal{A},B;D,G)\ \ \text{and}\
 \
 \widehat{X}^{\text{o}}_{\mathcal{A}}:=\widehat{\X}^{\text{o}}(\mathcal{A},B;D,G).
\end{equation*}}

\noindent{\it Proof of Step 2.}
First, we apply Proposition 2.5 to the domain $D$ which is locally
$\mathcal{C}^2$ smooth on an open  neighborhood of $P$ in $\partial D.$
Consequently, we may find an open neighborhood $U$ of $P$
satisfying (2.2) such that
 Proposition 2.5 is applicable there. In the sequel
     the notation $U,$ $U_1,$ $U_3,$ $\pi^{\C},$ $\pi,$ $V$ and $V_Q$ have the same
     meanings as in Proposition 2.5.
Now we can fix a ball $\mathcal{A} $ of $A$
  \begin{equation}
 \mathcal{A}:=A\cap B(P,r),
\end{equation}
where the radius $r$ is sufficiently small such that
$2\mathcal{A}\Subset A,$
$2\mathcal{A}\Subset \partial (U\cap D)\cap\partial
 V,$ $2\mathcal{A}\Subset U_3$  and  $\pi^{\C}(2\mathcal{A})\Subset U_1.$

For any $\delta$ small enough, by Proposition 3.5 we may find an open subset
$T_{\delta}$ of $U\cap D$ such that
\begin{equation}
\begin{split}
\omega(z,\mathcal{A},D)-\delta &\leq \omega(z,T_{\delta},D)
\leq \omega(z,\mathcal{A},D),\qquad \ z\in D,\\
\sup\limits_{T_{\delta}}\dist(\cdot,\mathcal{A})&<\delta\qquad\text{and}\qquad \pi^{\C}(T_{\delta})\Subset U_1.
\end{split}
\end{equation}

 A geometric argument based  on Proposition 2.5 and definition (5.2) shows that
 one may find $\delta_0>0$ small enough such that
for any $z\in D\cup A$ with $\dist(z,2\mathcal{A})<\delta_0,$ $z\in U_3$ and
 there is a unique
$Q_z\in U_1$ such that $z\in V_{Q_z}.$  In addition
by Part 4) of Proposition 2.5 we have \begin{equation}
   \dist\left(z,\partial V_{Q_z}\cap \frac{3}{2}\mathcal{A}\right)\leq   C_1\cdot
 \dist(z,\mathcal{A})
\end{equation}
for any  $z\in D\cup A$ with $\dist(z,\mathcal{A})<\delta_0.$

 On the other hand, combining
Corollary 2.4 and Propositions 2.5 and 3.1, we get the following estimate
\begin{equation}
\omega\left(z,\partial  V_{Q_z}\cap 2\mathcal{A}, V_{Q_z}\right)\leq
 C_2\cdot  \dist\left(z,\partial V_{Q_z}\cap\frac{3}{2} \mathcal{A}\right),
\end{equation}
where $C_1,\ C_2$ are finite constants independent of  $z\in D\cup A$ with
$\dist(z,\mathcal{A})<\delta_0.$

For each $Q\in \pi^{\C}(2\mathcal{A}),$ we apply Gonchar's Theorem (Theorem 2) to the function $f
\in \mathcal{C}\left(\X(\partial  V_{Q}\cap 2\mathcal{A},B; V_{Q},G)\right)
\cap  \mathcal{O}_s\left(\X^{\text{o}}(\partial  V_{Q}\cap 2\mathcal{A},B; V_{Q},G)\right) $
 in order to obtain an extension
function $\tilde{f}_Q \in \mathcal{C}\left(\widehat{\X}(\partial  V_{Q}\cap 2\mathcal{A}, B; V_{Q},G)\right)
\cap  \mathcal{O}\left(\widehat{\X}^{\text{o}}(\partial  V_{Q}\cap 2\mathcal{A}, B; V_{Q},G)\right) $
such that $\tilde{f}_Q=f$ on $\X(\partial  V_{Q}\cap 2\mathcal{A}, B; V_{Q},G).$

Collecting the family $(\tilde{f}_Q)_{Q\in \pi^{\C}(2\mathcal{A})},$
we obtain an extension function $\tilde{f}$ defined on the following set
\begin{multline}
\widetilde{X}_{\mathcal{A}}:=\left\lbrace
   (z,w)\in (D\cup 2\mathcal{A})\times (G\cup B):\ \exists Q\in \pi^{\C}(2\mathcal{A}),
\ z\in V_Q \right.\\
\left.\text{and}\ \omega(z,\partial  V_{Q_z}\cap 2\mathcal{A}, V_{Q_z})+
\omega(w,B, G)  <1     \right\rbrace,
\end{multline}
which is not necessarily open;  moreover
\begin{equation}
\tilde{f}=f\qquad\text{on}\ \X(\mathcal{A}, B;T_{\delta_0},G).
\end{equation}
In virtue of (5.3)--(5.5) we obtain an $\delta_0$ small enough such that
for $0<\delta<\delta_0,$
\begin{equation}
\omega(z,\partial  V_{Q_z}\cap 2\mathcal{A}, V_{Q_z})\leq C_1C_2\cdot\dist(z,\mathcal{A})<
C_1C_2\delta<1,\qquad z\in T_{\delta}.
\end{equation}
Therefore, by (5.6), (5.8) and Theorem 2   for $0<\delta<\delta_1:=\min\left\{\delta_0,\frac{1}{2C_1C_2}
\right\}$
 and $z\in T_{\delta},$
$\tilde{f}(z,\cdot)$ is holomorphic on the open set
\begin{equation}
G_{\delta}:=\left\lbrace w\in G:\ \omega(w,B,G)<1-2C_1C_2\delta \right\rbrace.
\end{equation}

We need the following
\begin{lem}
For any $(\zeta_0,w_0)\in \overline{\mathcal{A}}\times (G\cup B),$ there are an open
neighborhood $\mathcal{U}$ of $\zeta_0$ in $D\cup A$ and  an open
neighborhood $\mathcal{V}$ of $w_0$ in $G\cup B$ such that
$\mathcal{U}\times\mathcal{V}\subset \widetilde{X}_{\mathcal{A}}$ and
$\vert\tilde{f}\vert_{\mathcal{U}\times\mathcal{V}}<\infty.$
\end{lem}

\noindent{\it Proof of Lemma 5.2.}
Fix a point $(\zeta_0,w_0)\in \overline{\mathcal{A}}\times (G\cup B).$
 Let
\begin{equation}
\epsilon:=\frac{1-\omega(w_0,B,G)}{3}
\end{equation}
and choose an open neighborhood $\mathcal{V}$ of $w_0$ in $G\cup B$ such
that
\begin{equation}
\omega(w,B,G)<\omega(w_0,B,G)+\epsilon,\qquad  w\in\mathcal{V}.
\end{equation}
Moreover, choose  a sufficiently small open neighborhood $\mathcal{U}$ of $\zeta_0$ in $D\cup A$
 such that  by (5.4) and (5.5)  we have
\begin{equation}
 \omega(z,\partial  V_{Q_z}\cap 2\mathcal{A}, V_{Q_z}) <\epsilon,\qquad z\in \mathcal{U}.
\end{equation}

Next, by Proposition 3.3, we may find a   subdomain
$G^{\epsilon}$ of $G$ such that $w_0\in G^{\epsilon},$ $\overline{G^{\epsilon}}\subset G\cup B$
and $G^{\epsilon}$ is locally $\mathcal{C}^2$ smooth on the   open subset
$B_{\epsilon}:=\partial G\cap \partial G^{\epsilon}$ of $B,$
$B_{\epsilon}\Subset B$ and
\begin{equation}
\omega(w_0,B,G)\leq \omega(w_0,B_{\epsilon},G^{\epsilon})<\omega(w_0,B,G)+\epsilon.
\end{equation}
By shrinking $\mathcal{V}$ if necessary, we may assume that
\begin{equation}
  \omega(w,B_{\epsilon},G^{\epsilon})<\omega(w_0,B_{\epsilon},G^{\epsilon})+\epsilon,\qquad w\in\mathcal{V}.
\end{equation}
Since $f\in\mathcal{C}(X)$ and $2\mathcal{A}\Subset A,$ by shrinking $\mathcal{U},$ if necessary, we may
find a finite constant $M$ such that
\begin{equation}
\vert \tilde{f}\vert_{2\mathcal{A}\times G^{\epsilon}}<M
\qquad\text{and}\qquad
\vert \tilde{f}\vert_{\mathcal{U}\times B_{\epsilon}}<
M.
\end{equation}
Consequently, for each $Q\in \pi^{\C}(2\mathcal{A})$ we apply Gonchar's Theorem (Theorem 2) to the function $f
\in \mathcal{C}\left(\X(\partial  V_{Q}\cap 2\mathcal{A},B_{\epsilon}; V_{Q},G^{\epsilon})\right)
\cap  \mathcal{O}_s\left(\X^{\text{o}}(\partial  V_{Q}\cap 2\mathcal{A},B_{\epsilon}; V_{Q},G^{\epsilon})
\right) $
 and  obtain  the inequality  $\vert \tilde{f}\vert <M$  on the following set
\begin{multline*}
\widetilde{X}^{\epsilon}_{\mathcal{A}}:=\left\lbrace
   (z,w)\in (D\cup2\mathcal{A})\times (G^{\epsilon}\cup B_{\epsilon}):\ \exists Q\in \pi^{\C}(2\mathcal{A}),
\ z\in V_Q\right.\\
\left.\text{and}\ \omega(z,\partial  V_{Q_z}\cap 2\mathcal{A}, V_{Q_z})+
\omega(w,B_{\epsilon}, G^{\epsilon})  <1     \right\rbrace,
\end{multline*}
 On the other hand, using (5.6) and (5.11)--(5.14), we see that
 \begin{equation*}
\mathcal{U}\times\mathcal{V}\subset
\widetilde{X}^{\epsilon}_{\mathcal{A}}\subset \widetilde{X}_{\mathcal{A}}.
 \end{equation*}
Hence $\vert\tilde{f}\vert_{\mathcal{U}\times\mathcal{V}}<M,$ which
completes the proof of the lemma.
\hfill $\square$

\smallskip

Now for any
$0<\delta<\delta_1,$ we are able to apply Theorem 4.1 to the function
\begin{equation*}
\tilde{f}\in  \mathcal{C}_s\left( \X( T_{\delta},B; D,G_{\delta}
       )\right)\cap
 \mathcal{O}_s\left( \X^{\text{o}}( T_{\delta},B; D,G_{\delta}
       )\right)
\end{equation*}
 and obtain
a function $\hat{f}_{\delta}\in \mathcal{C}\left(\widehat{\X}( T_{\delta},B; D,G_{\delta}
       )\right)\cap
 \mathcal{O}\left(\widehat{\X}^{\text{o}}( T_{\delta},B; D,G_{\delta}
       )\right)
$ such that
\begin{equation}
\hat{f}_{\delta}=\tilde{f}\qquad \text{on}\ \X( T_{\delta},B; D,G_{\delta} ).
\end{equation}

We are now in a position to define the desired extension  function $\hat{f}.$
Indeed,
one  glues
$\left(\hat{f}_{\delta}\right)_{0<\delta<\delta_1}$ together to obtain
$\hat{f}$ in the following way
\begin{equation}
\hat{f}:=\lim\limits_{\delta\to 0} \hat{f}_{\delta}\qquad \text{on}\
\widehat{\X}\left(\mathcal{A},B;D,G\right) \setminus (\mathcal{A}\times (G\cup B)).
\end{equation}
One now checks that the limit (5.17) exists and possesses all the required
properties. This is an immediate consequence of the following
\begin{lem}
 For  any point $(z,w)\in \widehat{\X}\left(\mathcal{A},B;D,G\right)
 \setminus (\mathcal{A}\times (G\cup B))$ let $\delta_{z,w}$ be the unique positive
 number $\delta$  which verifies
 \begin{equation}
\delta+\frac{2C_1C_2\delta}{1-\delta}=1-\omega(z,\mathcal{A},D)-\omega(w,B,G).
\end{equation}
Then $f(z,w)=\hat{f}_{\delta}(z,w)$ for all $0<\delta<\delta_{z,w}.$
\end{lem}

\smallskip

\noindent{\it Proof of Lemma 5.3.}
Fix a point  $(z_0,w_0) \in \widehat{\X}\left(\mathcal{A},B;D,G\right)
  \setminus (\mathcal{A}\times (G\cup B)).$
Then by (5.3), (5.18) and the second identity of Part 1) of Proposition 3.8,
 for all $0<\delta <\delta_{z_0,w_0},$
\begin{equation*}
(z_0,w_0)\in \widehat{\X}( T_{\delta},B; D,G_{\delta}       ).
\end{equation*}
Therefore, for any $0<\delta^{'}<\delta<\delta_{z_0,w_0},$ $(z_0,w_0)$ is
contained in the following set
\begin{equation}
 \widehat{\X}( T_{\delta^{'}},B; D,G_{\delta^{'}}       )\cap
\widehat{\X}( T_{\delta},B; D,G_{\delta}       ).
\end{equation}
Let $H_{\delta}$ be the connected component of $G_{\delta}$ containing $w_0.$
Let $B_{\delta}$ be the largest open subset of $B$ such that
 $H_{\delta}$ is locally $\mathcal{C}^2$
smooth on $B_{\delta}.$
By  Part 2) of Proposition 3.8, $B_{\delta}$ is nonempty and $\omega(w,B_{\delta},H_{\delta})
=\omega(w,B,G_{\delta}),$  $w\in H_{\delta}.$   On the other hand, using (5.3)
and  the inclusion  $G_{\delta}\subset G_{\delta^{'}}    , $ we see that   the set  (5.19)
contains the set
\begin{equation}
\left\lbrace (z,w)\in D\times H_{\delta}:\ \omega(z,T_{\delta},D)+\omega(w, B_{\delta}
,H_{\delta})<1-\delta \right\rbrace.
\end{equation}
By Proposition 3.7 this open set is connected. Moreover, it contains the point $(z_0,w_0)$.
In addition by (5.7) and (5.16), one gets
\begin{equation*}
\hat{f}_{\delta^{'}}=\hat{f}_{\delta}=\tilde{f}=f \qquad \text{on}\ T_{\delta}\times
 B_{\delta}.
\end{equation*}
Applying Theorem 2.1, we deduce that  $\hat{f}_{\delta^{'}}=\hat{f}_{\delta}$
on the domain given by (5.20). In particular   $\hat{f}_{\delta^{'}}(z_0,w_0)=\hat{f}_{\delta}(z_0,w_0).$
This completes the proof.
\hfill $\square$

\smallskip

Another consequence of Lemma 5.3 is that
$\hat{f}\in\mathcal{O}\left(\widehat{\X}^{\text{o}}  \left(\mathcal{A},B;D,G\right)\right).$
Now we define $\hat{f}$ on $\mathcal{A}\times (G\cup B)$ as follows
\begin{equation}
\hat{f}:=f\qquad \text{on}\ \mathcal{A}\times (G\cup B).
\end{equation}
Thus $\hat{f}$ is well-defined on the whole
$ \widehat{\X}\left(\mathcal{A},B;D,G\right).$

To complete Step 2, it remains to show that  $\hat{f} \in
\mathcal{C}(\widehat{\X} \left(\mathcal{A},B;D,G\right))$  and $\hat{f}=f$ on
$\X \left(\mathcal{A},B;D,G\right).$

First we prove that $\hat{f}$ is continuous on $D\times B.$
For this let $(z_0,\eta_0)\in D\times B.$ By Proposition 3.2 there are an open
neighborhood $\mathcal{U}$ of $z_0$ in $D$ and an open neighborhood $\mathcal{V}$ of $\zeta_0$
in $G\cup B$  such that
\begin{equation*}
\lambda:=\sup\limits_{z\in \mathcal{U},\ w\in
\mathcal{V}}\left(\omega(z,\mathcal{A},D)+\omega(w,B,G)\right)
<1.
\end{equation*}
Now let $\delta>0$ be a positive number which verifies
 $\delta+\frac{2C_1C_2\delta}{1-\delta}<1- \lambda.$
Then  Lemma 5.3 implies that $f=\hat{f}_{\delta}$  on
$\mathcal{U}\times \mathcal{V}.$ Since by Theorem 4.1 we have  known that $\hat{f}_{\delta}$
is continuous on $D\times B,$
 so is $\hat{f}$
and moreover $\hat{f}=f$ on $D\times B.$
This, combined with (5.21) implies that  $\hat{f}=f$ on $\X \left(\mathcal{A},B;D,G\right).$

Finally, it remains to check the continuity of $\hat{f}$ on $\mathcal{A}\times (G\cup B).$
Fix a point $(\zeta_0,w_0)\in\mathcal{A}\times (G\cup B)$
and a number $0<\epsilon <1.$ From the hypothesis
$f\in\mathcal{C}(X)$ and by Lemma 5.2 it follows that there are an open connected neighborhood $\mathcal{U}$ of $\zeta_0$ in
$D\cup\mathcal{A},$  an open connected neighborhood $\mathcal{V}$ of $w_0$ in $G\cup B$
and a finite constant $M$ such that
\begin{equation}
\begin{split}
\left\vert f(\zeta_0,w_0)-f(\zeta, w)\right\vert&<\epsilon^{2},\qquad
 \zeta\in \mathcal{A}\cap\mathcal{U},\ w\in \mathcal{V},\\
 \vert \tilde{f}\vert_{\mathcal{U}\times\mathcal{V}}&< \frac{M}{2}.
\end{split}\end{equation}
Moreover, by shrinking $\mathcal{U}$ and $\mathcal{V},$ if necessary, and
applying Proposition 3.2, we may suppose that
\begin{equation*}
\sup\limits_{\mathcal{U}\times\mathcal{V}}\left(\omega(z,\mathcal{A},D)+\omega(w,B,G)
\right )<1.
\end{equation*}
Therefore, $\mathcal{U}\times\mathcal{V}\subset\widehat{\X}
\left(\mathcal{A},B;D,G\right).$
Moreover, by Lemma 5.3 and (5.17), there is an $\delta>0$ such that
$\hat{f}=\hat{f}_{\delta}=\tilde{f}$ on the nonempty open set
$(T_{\delta}\cap\mathcal{U})\times\mathcal{V}.$
Thus
\begin{equation}
\hat{f}=\tilde{f} \qquad \text{on} \ \mathcal{U}\times\mathcal{V}.
\end{equation}
By shrinking $\mathcal{U},$ if necessary, we may suppose that for all
$z\in\mathcal{U},$ there is exactly one point $\zeta_z\in\mathcal{A}$ such
that $\pi(\zeta_z)=\pi(z).$ By Part 4) of Proposition 2.5  we have
\begin{equation*}
z\in V_{Q_z}\qquad \zeta_z\in \mathcal{A}\cap \partial V_{Q_z}\qquad
\text{and}\qquad
\dist(z,\zeta_z)\approx \dist(z,\partial V_{Q_z}).
\end{equation*}
 Therefore, we are able to apply the Two-Constant Theorem   to the function
$\tilde{f}(\cdot,w)- \tilde{f}(\zeta_z,w)\in\mathcal{C}\left( (\partial V_{Q_z}\cap\mathcal{U})\cup
( V_{Q_z}\cap\mathcal{U})\right)
 \cap  \mathcal{O}( V_{Q_z}\cap\mathcal{U}),$
which is, by (5.22), bounded by $M$ for any $z\in \mathcal{U},$   $w\in \mathcal{V}.$

Consequently,   taking  (5.22) and (5.23) into account, we deduce that
\begin{equation*}
\vert \hat{f}(z,w)- \hat{f}(\zeta_z,w)\vert <\epsilon^{2(
1-\omega(z,\partial V_{Q_z}\cap\mathcal{U}, V_{Q_z}\cap\mathcal{U}))}
M^{\omega(z,\partial V_{Q_z}\cap\mathcal{U}, V_{Q_z}\cap\mathcal{U})}
\end{equation*}
for all $(z,w)\in\mathcal{U}\times\mathcal{V}.$ Thus for $(z,w)\in (D\cup\mathcal{A})\times (G\cup B)$
 sufficiently close to
$(\zeta_0,w_0)$ we have by Proposition 3.2 and (5.22),
\begin{equation*}
\vert \hat{f}(z,w)- \hat{f}(\zeta_0,w_0)\vert
\leq \vert \hat{f}(z,w)- \hat{f}(\zeta_z,w)\vert+
\vert \hat{f}(\zeta_z,w)- \hat{f}(\zeta_0,w_0)\vert<\frac{\epsilon}{2}+\frac{\epsilon}{2}= \epsilon,
\end{equation*}
which proves
the continuity of $\hat{f}$ at $(\zeta_0,w_0).$

Hence Step 2 is finished.
\hfill $\square$

\smallskip

\noindent{\bf Step 3:} {\it    The case where $G$ is a Jordan planar simply connected domain.
 }

\smallskip

 \noindent{\it Proof of Step 3.}
Fix a sequence $(A_k)_{k=1}^{\infty}$ of open subsets of $A$ such that
$A_k\Subset A_{k+1}$ and  $A_k\nearrow A$ as $k\nearrow\infty.$ By
Proposition 3.4, we obtain  $\widehat{\X}(A_k,B;
D,G)\nearrow\widehat{\X}(A,B; D,G).$   Using a routine uniqueness argument (Theorem 2.1 and Proposition 3.6)
 and the gluing procedure, we are reduced to the proof that
 for any $k,$ there is a function
 $\hat{f}_k\in
\mathcal{C}\left(\widehat{\X}(A_k,B; D,G)\right)\cap\mathcal{O}\left(
\widehat{\X}^{\text{o}}( A_k,B; D,G       )\right)$  satisfying
$\hat{f}_k=f$ on $\X( A_k,B; D,G       ).$

Now fix an $k\in \N.$ First we shall show that      one may  find
a sufficiently small positive number $\delta_0$ with the following
properties:

For any $0<\delta<\delta_0,$ there is a finite number of
 open balls  $\left( \mathcal{A}_j\right)^{N}_{j=1}$ of $A$
    with radius $\delta$ ($N$ depending on $\delta$) such that
\begin{itemize}
\item[(i)] $A_k\subset \bigcup_{j=1}^N  \mathcal{A}_j;$
\item[(ii)]  $\bigcup_{j=1}^N  2\mathcal{A}_j\subset A_{k+1},$ where
 $2\mathcal{A}_j$ is the ball with the same center as  $\mathcal{A}_j$
 but with double radius;
\item[(iii)]   for each $1\leq j\leq N,$ Step 2 can apply to the open ball
 $ 2 \mathcal{A}_j ;$ more precisely,
     Step 2 provides a function
 $\hat{f}_j\in
\mathcal{C}\left(\widehat{X}_{2\mathcal{A}_j}\right)\cap
\mathcal{O}\left( \widehat{X}^{\text{o}}_{2\mathcal{A }_j}\right)$  satisfying
$\hat{f}_j=f$ on $X_{2\mathcal{A}_j};$

\item[(iv)] for any $0<\delta<\delta_0$
there is an open subset
$T_{\delta}$ of $ D$ such that
\begin{eqnarray*}
\omega(z,A_k,D)-\delta &\leq& \omega(z,T_{\delta},D)
\leq \omega(z,A_k,D),\qquad \ z\in D,\\
\sup\limits_{T_{\delta}}\dist(\cdot,A_k)&<& r,
\end{eqnarray*}
for some $0<r:=r_{\delta}<\delta;$

\item[(v)] for any $0<\delta<\delta_0$
and $z\in
 T_{\delta}$ there is a unique nearest point
 $\zeta_z\in \bigcup_{j=1}^N\mathcal{A}_j$ such that
 $\dist(z,\zeta_z)=\dist(z,\partial D)$ and
 for any $1\leq j\leq N$ such that $\zeta_z\in \mathcal{A}_j$ we have
 \begin{equation*}
\sup\limits_{t\in[z,\zeta_z]} \omega(t,2\mathcal{A}_j,D)< \delta,
\end{equation*}
where $ [z,\zeta_z]$ denotes the real segment connecting $z$ to $\zeta_z.$
  \end{itemize}

Indeed, using the result of Step 2 and by a compactness argument
we see that one may find $\delta_0>0$  sufficiently small such that the properties (i)--(iii)
are fulfilled.

On the other hand, using Proposition 3.2 we see that there is an
$r:=r_{\delta}$ sufficiently small such that
\begin{equation}
\omega(z,2\mathcal{A}_j,D)< \delta,
\end{equation}
for all $1\leq j\leq N$ and $z\in D$ with $\dist(z,\mathcal{A}_j)<r.$

By examining carefully the proof of Proposition 3.5, we may arrange
$T_{\delta}$ in such a way  that property (iv) is fulfilled with $r$ given above.
It is also clear that when $r$ is sufficiently small,  the first assertion
of (v) is satisfied. Moreover the second one is an immediate consequence
of (5.24).

Thus we have shown that  all properties (i)--(v) are verified.

Next fix $0<\delta<\delta_0$ and put
\begin{equation}
G_{\delta}:=\left\lbrace w\in G:\   \omega(w,B,G)<1-2\delta \right\rbrace.
\end{equation}
We define a new function $\tilde{f}$ on $(T_{\delta}\cup A_k)\times
G_{\delta}$ as follows. For any $(z,w)\in (T_{\delta}\cup A_k)\times G_{\delta}$ let
\begin{equation}
\tilde{f}(z,w):=\hat{f}_j(z,w)
\end{equation}
for any $1\leq j\leq N$  such that $\zeta_z\in \mathcal{A}_j$ (see the
notation in (v) above).

First one checks that $\tilde{f}$ is well-defined. Indeed,
in virtue of (iv)--(v) and (5.25), for any $(z,w)\in (T_{\delta}\cup A_k)\times G_{\delta}$
there is at least an $j$ such that
$\zeta_z\in\mathcal{A}_j$ and  $(t,w)\in \widehat{X}_{2\mathcal{A}_j},$
for $t\in[z,\zeta_z].$

On the other hand, suppose that there is another
index $l$ such that $ \zeta_z\in\mathcal{A}_l.$
Observe that $\hat{f}_j=\hat{f}_l=f$ on $(\mathcal{A}_j\cap\mathcal{A}_l)\times G.$
Therefore, we may apply Theorem 2.1 and Proposition 3.6 and conclude that
$\hat{f}_j=\hat{f}_l$ on the connected component of $ \widehat{X}_{2\mathcal{A
}_j}\cap \widehat{X}_{2\mathcal{A}_l} $
which is locally $\mathcal{C}^2$ smooth on $ (\mathcal{A}_j\cap \mathcal{A}_l)\times G_{\delta}.$
However we have already shown in the previous paragraph  that $(t,w)\in  \widehat{X}_{2\mathcal{A
}_j}\cap \widehat{X}_{2\mathcal{A}_l}$ for $t\in[z,\zeta_z]$ and clearly
$(\zeta_z,w)\in (\mathcal{A}_j\cap\mathcal{A}_l)\times G.$
Consequently, the above mentioned  connected component contains the point
$(z,w).$
Thus $\hat{f}_j(z,w)=\hat{f}_l(z,w),$ and hence  the function $\tilde{f}$ is
 well-defined.

In virtue of (5.26), it is also clear that
\begin{equation*}
\tilde{f}\in \mathcal{C}\left((T_{\delta}\cup A_k)\times G_{\delta}\right)
 \cap\mathcal{O}\left(T_{\delta}\times G_{\delta}\right).
\end{equation*}
Let $\tilde{f}_{\delta}$ be the trace of $\tilde{f}$ on $\X(T_{\delta},B; D,G_{\delta}).$
Applying Theorem 4.1 to the function
$\tilde{f}_{\delta} \in \mathcal{C}\left(\X(T_{\delta},B; D,G_{\delta})\right)\cap
\mathcal{O}\left(\X^{\text{o}}(T_{\delta},B; D,G_{\delta}       )\right),$
 we obtain  an extension function
\begin{equation*}
\hat{f}_{\delta}\in
\mathcal{C}\left(\widehat{\X}(T_{\delta},B; D,G_{\delta})\right)\cap
\mathcal{O}\left(\widehat{\X}^{\text{o}}( T_{\delta},B; D,G_{\delta}       )\right)
\end{equation*}
 satisfying
$\hat{f}_{\delta}=f$ on $D\times B.$

Finally, one proceeds as in the end of Step 2. Observe that Lemma 5.3 is still valid
in the present context.  As in formula (5.17), one may glue
$\left(\hat{f}_{\delta}\right)_{0<\delta<\delta_0}$  to obtain an
extension function $\hat{f}:=\lim\limits_{\delta\to 0}\hat{f}_{\delta}$ which is holomorphic on
$\widehat{X}^{\text{o}}$ and continuous on $D\times B.$

Since $\tilde{f}\in \mathcal{C}\left((T_{\delta}\cup A_k)\times G_{\delta})\right)$
for $0<\delta<\delta_0,$
Lemma 5.3 in the present context also gives that
 $\hat{f} \in
\mathcal{C}\left(\widehat{\X}(A,B;D,G)\right).$

Hence Step 3 is complete.
\hfill $\square$

\noindent{\bf Step 4:} {\it We keep the hypothesis of Theorem 5.1  and
 prove the following local version of this theorem:}

\smallskip

{\it For any point $P\in B,$ there is  a ball
$\mathcal{B}$ in $B$ with center  $P$ such that
 the following property holds:  For any   function $f\in\mathcal{C}(X)\cap\mathcal{O}_s(X^{\text{o}}),$ there
  is a unique function
$\hat{f}\in\mathcal{C}(\widehat{X}_{\mathcal{B}})\cap \mathcal{O}(\widehat{X}^{\text{o}}_{\mathcal{B}})$
such that $\hat{f}=f$ on $X_{\mathcal{B}},$
where
 \begin{equation*}
 X_{\mathcal{B}}:=\X(A,\mathcal{B}; D,G), \
 \widehat{X}_{\mathcal{B}}:=\widehat{\X}(A,\mathcal{B};D,G)\ \ \text{and}\
 \
 \widehat{X}^{\text{o}}_{\mathcal{B}}:=\widehat{\X}^{\text{o}}(A,\mathcal{B};D,G).
\end{equation*}
}

 \noindent{\it Proof of Step 4.}
 We proceed using Step 3 in exactly the same way as we proved Step 2 using
 Theorem 2. Therefore we shall only indicate briefly the outline  of the
 proof.

First we apply Proposition 2.5 to the domain $G$ which is locally
$\mathcal{C}^2$ smooth on an open  neighborhood of $P$ in $\partial G.$
Consequently, we may find an open neighborhood $U$ of $P$
satisfying (2.2) such that
 Proposition 2.5 is applicable there. In the sequel
     the notation $U,$ $U_1,$ $\pi^{\C},$  $V$ and $V_Q$ have the same
     meanings as in Proposition 2.5.
Now we can fix a ball $\mathcal{B} $ of $B:$
  $\mathcal{B}:=B\cap B(P,r),$ where the radius $r$ is sufficiently small such that
$2\mathcal{B}\Subset B$ etc.

Arguing as in (5.2)--(5.3) we can choose an $\delta_0>0$ sufficiently small such that
for any $0<\delta<\delta_0$  there is an open subset
$S_{\delta}$ of $G$ satisfying
\begin{equation}
\begin{split}
\omega(w,\mathcal{B},G)-\delta &\leq \omega(w,S_{\delta},G)
\leq \omega(w,\mathcal{B},G),\qquad \ w\in G,\\
\sup\limits_{S_{\delta}}\dist(\cdot,\mathcal{B})&<\delta
 .
\end{split}
\end{equation}
 Arguing as in (5.4)--(5.8),
  there is a finite constant $C_3$ such that
\begin{equation}
\omega(w,\partial V_{Q_w}\cap 2\mathcal{B},V_{Q_w}) \leq C_3\dist(w,\mathcal{B})
\leq C_3\delta,
\end{equation}
for $0<\delta<\delta_0,$  $w\in S_{\delta}$ and $Q_w:=\pi^{\C}(w).$

Lemma 5.2 is still valid in the present context
making the obviously necessary changes in notation. There is only one important difference
between Step 2 and the present step. In Step 2 we apply Gonchar's Theorem
to (5.15) but in this step we appeal to Theorem 4.3.

Lemma 5.3 is also valid in the present context making the obviously necessary changes in
notation.

For each $Q\in \pi^{\C}(2\mathcal{B}),$ we apply the result of Step 3 to the function $f
\in \mathcal{C}\left(\X(A, \partial  V_{Q}\cap 2\mathcal{B};D, V_{Q})\right)
\cap  \mathcal{O}_s\left(\X^{\text{o}}(A,\partial  V_{Q}\cap 2\mathcal{B};D, V_{Q})\right) $ in order to obtain an extension
function $\hat{f}_Q\in \mathcal{C}\left(\widehat{\X}(A, \partial  V_{Q}\cap 2\mathcal{B};D, V_{Q})\right)
\cap  \mathcal{O}\left(\widehat{\X}^{\text{o}}(A,\partial  V_{Q}\cap 2\mathcal{B};D, V_{Q})\right) $
such that $\hat{f}_Q=f$ on $\X(A, \partial  V_{Q}\cap 2\mathcal{B};D, V_{Q}).$

Gluing the family $(\hat{f}_Q)_{Q\in \pi^{\C}(2\mathcal{B})},$
we obtain an extension function $\tilde{f}$ defined on
\begin{equation}
\left\lbrace   (z,w)\in D\times G:\ \exists Q\in \pi^{\C}(2\mathcal{B}),
\ w\in V_Q\ \text{and}\ \omega(z,A,D)+\omega (w,\partial  V_{Q}\cap 2\mathcal{B}, V_{Q})
   <1     \right\rbrace.
\end{equation}
For $0<\delta<\delta_0$  put
\begin{equation}
D_{\delta}:=\left\lbrace \omega(z,A,D)<1-2C_3\delta\right\rbrace .
\end{equation}

As in Step 2, taking (5.27)--(5.30) into account we see that
\begin{equation*}
\tilde{f}\in \mathcal{C}_s\left(\X(A, S_{\delta};D_{\delta}, G)\right)
\cap  \mathcal{O}_s\left(\X^{\text{o}}(A, S_{\delta};D_{\delta}, G
)\right).
\end{equation*}
Therefore, we are in a position to apply Theorem 4.1 and obtain an extension function
\begin{equation*}
\hat{f}_{\delta}\in \mathcal{C}\left(\widehat{\X}(A, S_{\delta};D_{\delta}, G)\right)
\cap  \mathcal{O}\left(\widehat{\X}^{\text{o}}(A, S_{\delta};D_{\delta}, G
)\right).
\end{equation*}

Using (5.17) we may glue $\left(\hat{f}_{\delta}\right)_{0<\delta<\delta_0}$ together
in order to obtain the desired extension function $\hat{f}.$
The rest of the proof follows along the same lines as in Step 2 making use
of Two-Constant Theorem and Lemmas 5.2 and 5.3.
This finishes Step 4.
    \hfill $\square$

\smallskip

\noindent{\bf Step 5:} {\it  The general case.}

\smallskip

 The same argument which has been
 used to go from
 Step 2 to Step 3  will enable us  to go from Step 4 to Step 5.
 Consequently,  there is an extension function
 $\hat{f}\in\mathcal{C}(\widehat{X})\cap\mathcal{O}_s(\widehat{X}^{\text{o}})$
 such that $\hat{f}=f$ on $X.$ It is also clear that $\hat{f}$ is uniquely
 determined.
 Finally, it remains to establish estimate (5.1).  But it follows immediately from
 Theorem 4.3.

This completes the last step of the proof.
\end{proof}
\section{Proof of the Main Theorem and concluding remarks}
In order to prove the Main Theorem,
we proceed by induction (I) on $N\geq 2.$ Suppose the Main Theorem is true
for $N-1\geq 2.$ We have to discuss the case of an $N$-fold cross
$X:=\X(A_1,\ldots,A_N;D_1,\ldots,D_N),$ where $D_1,\ldots,D_N$ are
pseudoconvex domains and $A_1,\ldots,A_N$ are open subsets of $
\partial D_1,\ldots,\partial D_N$ such that $D_j$ is  locally $\mathcal{C}^2$
smooth  on $A_j$ $(1\leq j\leq N).$
Fix an  $f\in\mathcal{C}(X)\cap\mathcal{O}_s(X^{\text{o}}).$

We proceed again by induction (II) on the positive integer $j$ $(1\leq j\leq N)$ such that
$D_j,\ldots,D_N$ are Jordan planar domains.

For $j=1,$ we are reduced to Theorem 2.

 Suppose the Main Theorem is true
for the case where  $D_{j-1},\ldots,D_N$ are Jordan planar domains $(j\geq 2)$.
 We have to discuss the case where $D_j,\ldots,D_N$ are Jordan
planar domains. The  proof given below follows essentially the schema of
that of Theorem 5.1. It is divided into three steps.

\smallskip

\noindent{\bf Step 1:} {\it Reduction to the case where $D_1,\ldots,D_{j-1}$  are
bounded pseudoconvex domains.
}

\smallskip

 \noindent{\it Proof of Step 1.}  We proceeds in exactly the same way as
 in Step 1 of Theorem 5.1.
  This completes Step 1.
\hfill $\square$

\smallskip

{\it From now on we assume that the hypothesis of Step 1 is fulfilled. }

\smallskip

\noindent{\bf Step 2:} {\it   We  prove the following local version of the Main Theorem:}

\smallskip

{\it For any point $P\in A_1,$ there is  a ball
$\mathcal{A}$ in $A_1$ with center  $P$ such that
 the following property holds:  For any   function $f\in\mathcal{C}(X)\cap\mathcal{O}_s(X^{\text{o}}),$ there
  is a unique function
$\hat{f}\in\mathcal{C}(\widehat{X}_{\mathcal{A}})\cap \mathcal{O}(\widehat{X}^{\text{o}}_{\mathcal{A}})$
such that $\hat{f}=f$ on $X_{\mathcal{A}},$
where
 \begin{eqnarray*}
 X_{\mathcal{A}}&:=&\X\left(\mathcal{A},A_2,\ldots,A_N; D_1,\ldots,D_N\right),
 \qquad
 \widehat{X}_{\mathcal{A}}:=\widehat{\X}\left(\mathcal{A},A_2,\ldots,A_N; D_1,\ldots,D_N\right), \\
 \widehat{X}^{\text{o}}_{\mathcal{A}}&:=&\widehat{\X}^{\text{o}}\left(\mathcal{A},A_2,\ldots,A_N; D_1,\ldots,D_N\right) .
\end{eqnarray*}
}

 \noindent{\it Proof of Step 2.}  As in Step 2 in the proof  of Theorem 5.1  we first
 apply Proposition 2.5 to the domain $D_1$ which is locally
$\mathcal{C}^2$ smooth on an open  neighborhood of $P$ in $\partial D_1.$
Consequently, we may find an open neighborhood $U$ of $P$
satisfying (2.2) such that
 Proposition 2.5 is applicable there. In the sequel
     the notation $U,$ $U_1,$ $\pi^{\C},$  $V$ and $V_Q$ have the same
     meanings as in Proposition 2.5.
Now we can fix a ball $\mathcal{A} $ of $A_1:$ $
 \mathcal{A}:=A_1\cap B(P,r),$
where the radius $r$ is sufficiently small such that
$2\mathcal{A}\Subset A_1$
etc.

Arguing as in (5.2)--(5.3) and (5.4)--(5.8),
we can choose an $\delta_0>0$ sufficiently small
and a finite constant $C$ such that for any $0<\delta<\delta_0$  there is an open subset
$T^1_{\delta}$ of $D_1$ satisfying
\begin{equation}
\begin{split}
\omega(z_1,\mathcal{A},D_1)-\delta &\leq \omega(z_1,T^1_{\delta},D_1)
\leq \omega(z_1,\mathcal{A},D_1),\qquad \ z_1\in D_1,\\
\sup\limits_{T^1_{\delta}}\dist(\cdot,\mathcal{A})&<\frac{\delta}{C}
 ,\qquad  T^{1}_{\delta}\subset
 \bigcup\limits_{Q\in\pi^{\C}(2\mathcal{A})}V_Q,
\end{split}
\end{equation}
 and
\begin{equation}
\omega(z_1,\partial V_{Q_{z_1}}\cap 2\mathcal{A},V_{Q_{z_1}}) \leq C\dist(z_1,\mathcal{A})
\leq \delta,
\end{equation}
for $0<\delta<\delta_0,$  $z_1\in T^1_{\delta}$ and $Q_{z_1}:=\pi^{\C}(z_1).$

Similarly, for each $2\leq k\leq N$
there is an open subset
$T^k_{\delta}$ of $D_k$ satisfying
\begin{equation}
\omega(z_k,A_k,D_k)-\delta \leq \omega(z_k,T^k_{\delta},D_k)
\leq \omega(z_k,A_k,D_k),\qquad \ z_k\in D_k.
\end{equation}

For each $Q\in \pi^{\C}(2\mathcal{A}),$ we apply the reduction  assumption (II) to the function
\begin{multline*}f
\in \mathcal{C}\left(\X\left(2\mathcal{A}\cap \partial  V_{Q},A_2,\ldots,A_N; V_{Q},D_2,\ldots,D_N\right)\right)
\\
\cap  \mathcal{O}_s\left(\X^{\text{o}}\left(2\mathcal{A}\cap \partial  V_{Q},A_2,\ldots,A_N;
 V_{Q},D_2,\ldots,D_N\right)\right)
\end{multline*}
 in order to obtain an extension
function \begin{multline}
\hat{f}_Q\in \mathcal{C}\left(\widehat{\X}\left(2\mathcal{A}\cap \partial  V_{Q},A_2,\ldots,A_N; V_{Q},D_2,\ldots,D_N\right)\right)
 \\\cap
 \mathcal{O}\left(\widehat{\X}^{\text{o}}\left(2\mathcal{A}\cap \partial  V_{Q},A_2,\ldots,A_N;
 V_{Q},D_2,\ldots,D_N\right)\right)
\end{multline}
such that
\begin{equation}
\hat{f}_Q=f\qquad\text{on}\  \X(\mathcal{A}\cap \partial  V_{Q},A_2,\ldots,A_N; V_{Q},D_2,\ldots,D_N).
\end{equation}
Collecting  the family $(\hat{f}_Q)_{Q\in \pi^{\C}(2\mathcal{A})},$
we obtain an extension function $\tilde{f}$ defined on
\begin{multline}
\left\lbrace   (z_1,\ldots,z_N)\in D_1\times\cdots\times D_N:\ \exists Q\in \pi^{\C}(2\mathcal{A}),
\ z_1\in V_Q\ \text{and}\ \right.\\
\left.\omega(z_1,2\mathcal{A}\cap \partial V_Q,V_Q)
+\sum\limits_{k=2}^N\omega (z_k, A_k, D_k)
   <1     \right\rbrace
\end{multline}
which satisfies
\begin{equation}
\tilde{f}=f \qquad\text{on}\ \X\left(\mathcal{A},A_2,\ldots,A_N;T^1_{\delta_0},D_2,\ldots,D_N\right).
\end{equation}

For $0\leq\delta<\delta_0$  put
\begin{multline}
D^{'}_{\delta}:=\left\lbrace  (z_2,\ldots,z_N)\in D_2\times\cdots\times D_N:
\ \omega(z_2,A_2,D_2)+\omega(z_2,T^{3}_{\delta},D_3)+\right.
\\\left.\ \cdots+\omega(z_N,T^N_{\delta},D_N)<1- N\delta
\right\rbrace
\end{multline}
and
\begin{equation}
D^k_{\delta}:=\left\lbrace z_k\in D_k:\
\omega(z_k,A_{k},D_k)<1-N\delta\right\rbrace ,\qquad 1\leq k\leq N.
\end{equation}
Consequently, in virtue of (6.1)--(6.4) and (6.6)
for any fixed
 $z_1\in T^{1}_{\delta}$ and $0<\delta <\delta_0,$
  the restriction $\tilde{f}(z_1,\cdots)$
is holomorphic on  $ D^{'}_{\delta}.$

On the other hand, for any $a_2\in A_2,$ by the reduction assumption (I)
for an
$(N-1)$-fold cross, we obtain an extension $\hat{f}_{a_2}$  such that
\begin{equation}
\hat{f}_{a_2}\in\mathcal{C}\left(\widehat{\X}(A_1,A_3,\ldots,A_N;D_1,\ldots,D_N)\right)
\cap\mathcal{O}\left(\widehat{\X}^{\text{o}}(A_1,A_3,\ldots,A_N;D_1,\ldots,D_N)\right)
\end{equation}
and
\begin{multline}
\hat{f}_{a_2}(z_1,z_3,\ldots,z_N)=f(z_1,a_2,z_3,\ldots,z_N),\\
\qquad (z_1,z_3,\ldots,z_N)\in
\X(A_1,A_3,\ldots,A_N;D_1,D_3,\ldots,D_N).
\end{multline}

Observe that by (6.1)--(6.3), (6.6), and (6.9)--(6.10), for $0<\delta<\delta_0$ sufficiently
small,
the domain of definition of $\hat{f}_{a_2}$ $(a_2\in A_2)$ contains
$D^{1}_{\delta}\times  T^3_{\delta}\times\cdots\times
T^N_{\delta}$ and that of $\tilde{f}$ contains
$T^{1}_{\delta}\times A_2\times  T^3_{\delta}\times\cdots\times
T^N_{\delta}.$ Next
we would like to prove that for $0<\delta<\delta_0$ sufficiently small
and $a_2\in A_2,$
\begin{multline}
\tilde{f}(z_1,a_2,z_3,\ldots,z_N)=\hat{f}_{a_2}(z_1,z_3,\ldots,z_N),\\
(z_1,z_3,\ldots,z_N)\in T^{1}_{\delta}\times  T^3_{\delta}\times\cdots\times
T^N_{\delta}.
\end{multline}
Indeed, in virtue  (6.5) and (6.11) and by applying  the reduction assumption (I)
to $\hat{f}_{a_2}$ and the reduction assumption  (II) to $\hat{f}_Q$ for
 any $Q\in\pi^{\C}(2\mathcal{A})$ we know that
\begin{multline*}
\tilde{f}(z_1,a_2,z_3,\ldots,z_N)=\hat{f}_{Q}(z_1,a_2,z_3,\ldots,z_N)=
f(z_1,a_2,z_3,\ldots,z_N)=\hat{f}_{a_2}(z_1,z_3,\ldots,z_N),\\
\qquad z_1\in \mathcal{A}\cap\partial V_Q,\ a_2\in A_2\quad\text{and}\qquad
(z_3,\ldots,z_N)\in \X(A_3,\ldots,A_N;D_3,\ldots,D_N)
.
\end{multline*}
This proves (6.12). Consequently, we can define a new function
$\tilde{f}_{\delta}$ on $ \X\left(T^1_{\delta},A_2\times T^3_{\delta}\times\cdots\times
T^N_{\delta};D^1_{\delta}, D^{'}_{\delta}\right)$ as follows
\begin{equation}
  \tilde{f}_{\delta}:=
\begin{cases}
\tilde{f},
  & \qquad\text{on}\ T^{1}_{\delta}\times D^{'}_{\delta},\\
 \hat{f}_{a_2}, &   \qquad\text{on}\
 T^1_{\delta}\times\{a_2\}\times T^3_{\delta}\times\cdots\times
T^N_{\delta},\ a_2\in A_2         .
\end{cases}
\end{equation}
We need the following lemmas
\begin{lem} The following assertions hold:
\\
1) $\tilde{f}_{\delta}$ is locally bounded on
$ \X\left(T^1_{\delta},A_2\times T^3_{\delta}\times\cdots\times
T^N_{\delta};D^1_{\delta}, D^{'}_{\delta}\right);$\\
2) $\tilde{f}_{\delta}$ is locally bounded on
$T^1_{\delta}\times \left(D^{'}_{\delta}\cup \left( \overline{D^{'}_{\delta}}\cap
\X(A_2,\ldots,A_N;D_2,\ldots,D_N)       \right)       \right)$
and  $\tilde{f}_{\delta}(z_1,\cdots)\in \mathcal{C}\left(
D^{'}_{\delta}\cup \left( \overline{D^{'}_{\delta}}\cap
\X(A_2,\ldots,A_N;D_2,\ldots,D_N)
 \right )\right)$ for any $z_1\in T^1_{\delta}.$
\end{lem}

\noindent{\it Proof of Lemma 6.1.} It follows along the same lines as that of Lemma
5.2. Therefore, we only indicate a crucial difference.
In Lemma 5.2 we appeal to Gonchar's Theorem but in the present lemma
we apply the hypothesis of induction (II). This completes the proof.
\hfill  $\square$

\smallskip

\begin{lem} Let $D_2$ be a bounded open set and let
$A_2$ be an open set  of $\partial D_2$ such that $D_2$ is locally
$\mathcal{C}^2$ smooth on $A_2.$
Let $T^{k}\subset D_k\Subset \C^{n_k},$ $D_k$ a domain and $T^k$ locally
pluriregular, $k=3,\ldots,N,$ $N\geq 3.$
Put
\begin{multline*}
D^{'}:=\left\lbrace z^{'}=(z_2,\ldots,z_N)\in D_2\times\cdots\times D_N:\
\omega(z_2,A_2,D_2)\right.\\
\left. +\sum\limits_{k=3}^N\omega(z_k,T^k,D_k)<1
     \right\rbrace.
\end{multline*}
Then
\begin{equation*}
\omega\left(z^{'}, A_2\times T^3\times\cdots\times
T^N,D^{'}\right)=\omega(z_2,A_2,D_2)+\sum\limits_{k=3}^N\omega(z_k,T^k,D_k).
\end{equation*}
\end{lem}
\noindent{\it Proof of Lemma 6.2.} We argue as in the proof of Lemma 3(b)
in \cite{jp2} making use of Part 1) of Proposition 3.8.\hfill $\square$

\smallskip

We now come back to the proof of the Main Theorem. Applying Lemma 6.2
and Part 1) of  Proposition 3.8, we see that

\begin{multline}
\omega\left((z^{'},A_2\times
T^3_{\delta}\times \cdots\times T^N_{\delta},D^{'}_{\delta}\right)
=\frac{1}{1-N\delta}\left(
\omega(z_2,A_2,D_2) +\sum\limits_{k=3}^N\omega(z_k,T^k_{\delta},D_k)\right)
\end{multline}
 for any $z^{'}\in D^{'}_{\delta}.$

 To summarize what has been done so far: for any $0<\delta<\delta_0$ sufficiently small, we  obtain,
 by Part 1) of Lemma 6.1,  a function
$\tilde{f}_{\delta}$ defined on a mixed cross
\begin{equation*}
  \tilde{f}_{\delta}\in\mathcal{C}_s\left(     \X\left(T^1_{\delta},A_2\times T^3_{\delta}\times\cdots\times
T^N_{\delta};D^1_{\delta}, D^{'}_{\delta}\right)\right)
\cap \mathcal{O}_{s}\left(     \X^{\text{o}}\left(T^1_{\delta},A_2\times T^3_{\delta}\times\cdots\times
T^N_{\delta};D^1_{\delta}, D^{'}_{\delta}\right) \right) .
\end{equation*}
Applying Theorem 4.1 to $\tilde{f}_{\delta}$ we obtain an extension function  $\tilde{\tilde{f}}_{\delta}$
of  $\tilde{f}_{\delta}$
such that
\begin{gather}
 \tilde{\tilde{f}}_{\delta}\in\mathcal{C}\left(    \widehat{\X}\left(T^1_{\delta},A_2\times T^3_{\delta}\times\cdots\times
T^N_{\delta};D^1_{\delta}, D^{'}_{\delta}\right)\right)
\cap \mathcal{O}\left(     \widehat{\X}^{\text{o}}\left(T^1_{\delta},A_2\times T^3_{\delta}\times\cdots\times
T^N_{\delta};D^1_{\delta}, D^{'}_{\delta}\right) \right).
\end{gather}

 In virtue of (6.4)--(6.9) and (6.13)--(6.15) and by Part 2) of Lemma 6.1,
  we can apply Part 2) of Theorem 4.1
 and conclude that $\tilde{\tilde{f}}_{\delta}$ can be   continuously extended
 to a new function
\begin{equation}
\hat{f}_{\delta}\in\mathcal{C}(\widehat{X}_{\delta}),
\end{equation}
where
\begin{multline}
\widehat{X}_{\delta}:=
\widehat{\X}\left(T^1_{\delta},A_2\times T^3_{\delta}\times\cdots\times
T^N_{\delta};D^1_{\delta}, D^{'}_{\delta}\right)
\\
\cup\left(\overline{\widehat{\X}}\left(T^1_{\delta},A_2\times T^3_{\delta}\times\cdots\times
T^N_{\delta};D^1_{\delta}, D^{'}_{\delta}\right)\cap\left( D^{1}_{\delta}\times
\X\left(A_2,\ldots,A_N;D_2,\ldots,D_N \right)\right)
\right).
\end{multline}

We are now in a position to define the desired extension  function $\hat{f}.$
Indeed,
one  glues
$\left(\hat{f}_{\delta}\right)_{0<\delta<\delta_0}$ together to obtain
$\hat{f}$ in the following way
\begin{equation}
\hat{f}:=\lim\limits_{\delta\to 0} \hat{f}_{\delta}\qquad \text{on}\
\widehat{X}_{\mathcal{A}} \setminus \left(\mathcal{A}\times \X\left(A_2,\ldots,A_N;D_2,\ldots,D_N \right)\right).
\end{equation}
One now checks that the limit (6.18) exists and does posses all the required
properties. This is an immediate consequence of the following
\begin{lem}
 For  any point
 \begin{equation*}
 z=(z_1,\ldots,z_N)\in\widehat{X}_{\mathcal{A}} \setminus
  \left(\mathcal{A}\times \X\left(A_2,\ldots,A_N;D_2,\ldots,D_N \right)\right)
  \end{equation*}
  let $\delta_{z}$ be the unique positive number $\delta$  which verifies
 \begin{equation}
\delta+\frac{2N\delta}{1-N\delta}=1-\omega(z_1,\mathcal{A},D_1)-\sum\limits_{k=2}^N\omega(z_k,A_k,D_k).
\end{equation}
Then $f(z)=\hat{f}_{\delta}(z)$ for all $0<\delta<\delta_{z}.$
\end{lem}

\smallskip

\noindent{\it Proof of Lemma 6.3.}
Fix a point  $z^0=(z^0_1,\ldots,z^0_N) \in\widehat{X}_{\mathcal{A}} \setminus \left(\mathcal{A}\times
\X\left(A_2,\ldots,A_N;D_2,\ldots,D_N \right)
 \right).$
Then by (6.1)--(6.3) and (6.14)--(6.18) and the second identity of Part 1) of Proposition 3.8,
 for all $0<\delta <\delta_{z^0},$ we have
$
z^0\in \widehat{X}_{\delta}.
$ In particular,
 $\bigcup\limits_{0<\delta<\delta_0}\widehat{X}_{\delta}=\widehat{X}_{\mathcal{A}}
 \setminus\left (\mathcal{A}\times \X\left(A_2,\ldots,A_N;D_2,\ldots,D_N \right)\right).$

Therefore, for any $0<\delta^{'}<\delta<\delta_{z^0},$ $z^0$ is
contained in the following set
\begin{equation*}
 \widehat{X}_{\delta}\cap
\widehat{X}_{\delta^{'}}      .
\end{equation*}
Let $H_{\delta}$ be the connected component of $D^{'}_{\delta}$ containing $(z^0_2,\ldots,z^0_n).$
Let $B_{\delta}$ be the largest open subset of $B$ such that
 $H_{\delta}$ is locally $\mathcal{C}^2$
smooth on $B_{\delta}.$
By (6.1)--(6.3) and (6.19) and
arguing as in the proof of Lemma 5.3, we see that the above intersection contains
  the set
\begin{equation}
\left\lbrace z=(z_1,z^{'})\in D^{1}_{\delta}\times H_{\delta}:\ \omega(z_1,T^1_{\delta},D^{1}_{\delta})
+\omega(z^{'}, B_{\delta}
,H_{\delta})<1-\delta \right\rbrace.
\end{equation}
By Proposition 3.7 this open set is connected. Moreover, it contains the point $z^0$.
In addition we deduce  from  (6.17) that
\begin{equation*}
\hat{f}_{\delta^{'}}=\hat{f}_{\delta}=\tilde{f}=f \qquad \text{on}\ T^1_{\delta^{'}}\times
 B_{\delta^{'}}.
\end{equation*}
Applying Theorem 2.1, we deduce that  $\hat{f}_{\delta^{'}}=\hat{f}_{\delta}$
on the domain given by (6.20). In particular   $\hat{f}_{\delta^{'}}(z^0)=\hat{f}_{\delta}(z^0).$
This completes the proof.
\hfill $\square$

An immediate consequence of Lemma 6.3 is that
$\hat{f}\in\mathcal{O}\left(\widehat{X}^{\text{o}}_{\mathcal{A}}\right).$
Now we define $\hat{f}$ on $\mathcal{A}\times\X\left(A_2,\ldots,A_N;D_2,\ldots,D_N \right)
 $ as follows
\begin{equation*}
\hat{f}:=f\qquad \text{on}\ \mathcal{A}\times \X\left(A_2,\ldots,A_N;D_2,\ldots,D_N
\right).
\end{equation*}
Thus $\hat{f}$ is well-defined on the whole
$ \widehat{X}_{\mathcal{A}}$ and
\begin{equation}
\hat{f}\in\mathcal{C}\left(D_1\times  \X\left(A_2,\ldots,A_N;D_2,\ldots,D_N
\right)\right).
\end{equation}
To complete Step 2, it remains to show that  $\hat{f} \in
\mathcal{C}(\widehat{X}_{\mathcal{A}})$  and $\hat{f}=f$ on
$X_{\mathcal{A}}.$ For this purpose
we do the following trick.

We replace $D_1$ by $D_j$  $(j=2,\ldots,N)$  and proceed as above.
For example, if we replace $D_1$ by $D_2,$ then we obtain
a new extension function $\hat{\hat{f}}$  such that  in virtue of (6.21)
\begin{equation}
\hat{\hat{f}}\in\mathcal{C}\left(D_2\times  \X\left(A_1,A_3\ldots,A_N;D_1,D_3,\ldots,D_N
\right)\right).
\end{equation}
Next, using identities (6.12), (6.13), (6.16) and (6.18) and applying Theorem 2, we see that
the value of  $\hat{f}$ and $\hat{\hat{f}}$ can be uniquely determined
on  $T^1_{\delta}\times T^{2}_{\delta}\times  \cdots\times
T^N_{\delta}$
from the value of $f$ on $ \mathcal{A}\times A_2\times\cdots\times A_N $
for any sufficiently small $\delta>0.$  Thus
\begin{equation*}
\hat{\hat{f}}=\hat{f}\qquad\text{on}\
 T^1_{\delta}\times T^{2}_{\delta}\times  \cdots\times
T^N_{\delta}.
\end{equation*}
Hence $\hat{\hat{f}}=\hat{f}$
on $\widehat{X}_{\mathcal{A}}$ since  $\widehat{X}_{\mathcal{A}}$ is a
domain by Proposition 3.6.
Therefore, in virtue of (6.21), (6.22) and similar conclusions when $D_1$
is replaced by $D_3,\ldots,D_N,$
we conclude that
\begin{equation*}
\hat{f}\in \mathcal{C}\left(\widehat{X}_{\mathcal{A}}\setminus \left(\mathcal{A}\times A_2\times
\cdots\times A_N\right)\right).
\end{equation*}
Therefore,  Step 2 will be finished if we can  prove that
$\hat{f}$ is continuous on $ \mathcal{A}\times A_2\times
\cdots\times A_N.$

To do this fix a point $a=(a_1,\ldots,a_N)\in \mathcal{A}\times A_2\times
\cdots\times A_N$ and  an arbitrary number $\epsilon >0.$
 Next, we apply Proposition 2.5 to each domain $D_j$ which is locally
$\mathcal{C}^2$ smooth on an open  neighborhood of $a_j,$  $j=1,\ldots,N.$
Consequently, we may find an open neighborhood $U^{j}$ of $a_j$
satisfying (2.2) such that
 Proposition 2.5 is applicable there. In the sequel
     the notation $U^j,$ $U^j_1,$ $\pi^{\C,j},$  $V^j$ and $V^j_Q$ have the same
     meanings for $a_j$ as $U,\ U_1,\ \pi^{\C},\ V$ and $V_Q$ do for $P$ in Proposition 2.5.

Since $f\in\mathcal{C}(X),$ by shrinking $U^j,$ if necessary,  we may
assume without loss of generality that
\begin{equation}
\vert f(\zeta)-f(\eta)\vert <\frac{\epsilon}{2},\qquad
\zeta,\eta\in \X\left(\mathcal{A}\cap U^1, A_2\cap U^2,\ldots,A_N\cap U^N;
D_1\cap U^1,\ldots,D_N\cap U^N \right).
\end{equation}

Let $z=(z_1,\ldots,z_N)$ be an arbitrary point of $U^1\times\cdots\times U^N$
and put $Q_j:=\pi^{\C,j}(z_j).$
Then, in virtue of the hypothesis on $f,$ we may apply Theorem 2 to
 \begin{multline*}
f\in \mathcal{C}\left(\X\left(\mathcal{A}\cap \partial  V^1_{Q_1},A_2\cap \partial  V^2_{Q_2},
\ldots,A_N\cap \partial  V^N_{Q_N};
V^1_{Q_1},\ldots,V^N_{Q_N}\right)\right)\\
 \cap
 \mathcal{O}\left(\X^{\text{o}}\left(\mathcal{A}\cap \partial  V^1_{Q_1},A_2\cap \partial  V^2_{Q_2},
\ldots,A_N\cap \partial  V^N_{Q_N}; V^1_{Q_1},\ldots,V^N_{Q_N}\right)
 \right).
\end{multline*}
Consequently, taking  into account (6.23) and the above construction of the extension
function $\hat{f},$  we deduce that
\begin{multline*}
\vert \hat{f}(z)-f(\zeta)\vert <\frac{\epsilon}{2},
\qquad \zeta\in\X\left(\mathcal{A}\cap \partial  V^1_{Q_1},A_2\cap \partial  V^2_{Q_2},
\ldots,A_N\cap \partial  V^N_{Q_N}; V^1_{Q_1},\ldots,V^N_{Q_N}\right).
\end{multline*}
Hence, fixing any $\zeta$ as above and applying again (6.23), we get
\begin{equation*}
\vert \hat{f}(z)-f(a)\vert \leq\vert \hat{f}(z)-f(\zeta)\vert
+\vert f(\zeta)-f(a)\vert  <\frac{\epsilon}{2}+\frac{\epsilon}{2}
<\epsilon,
\end{equation*}
which proves the continuity of $\hat{f}$ at $a.$
Hence the remaining assertion  of Step 2 is proved.

Thus
 the proofs for the induction  (I) and (II)   are complete in this second step.
\hfill  $\square$

\smallskip

\noindent{\bf Step 3:} {\it  The general case.}

\smallskip

  The same argument which has been
 used to go from
 Step 2 to Step 3 in the proof of Theorem 5 will enable us  to go from Step 2 to Step 3
 in the present context.
Consequently,  there is an extension function
 $\hat{f}\in\mathcal{C}(\widehat{X})\cap\mathcal{O}_s(\widehat{X}^{\text{o}})$
 such that $\hat{f}=f$ on $X.$ It is also clear that $\hat{f}$ is uniquely
 determined.
 Finally, it remains to establish estimate (1.3). We have already proved the existence and uniqueness
 of the Main Theorem. Using  this result we argue as  in the proof of  Theorem
 4.2 in order to obtain (1.3).
This completes the last step of the proof.

 Hence the Main Theorem is proved.
\hfill $\square$

\medskip

Finally, we conclude this paper by some open remarks and open questions.\\

\smallskip

\noindent 1.  It seems to be of interest to establish the Main Theorem under weaker
assumptions than  the continuity of $f,$   the smoothness of $D_j$ on $A_j,$
 and  the regularity of the set $A_j$  $j=1,\ldots,N,$ etc. We postpone this issue to an ongoing work.\\

\noindent 2. Does the Main Theorem still hold if we only assume that
$A_j$ is of positive $(2n_j-1)$-Hausdorff measure,  $j=1,\ldots,N$?

\end{document}